\documentclass[a4paper,preprint,12pt]{elsarticle} 
\usepackage[english]{babel} 
\usepackage[latin1]{inputenc} 
\usepackage{amssymb} 
\usepackage{amsmath}    
\usepackage{pb-diagram}      
\usepackage{epstopdf}                    
\usepackage{graphicx}           
 \usepackage{epsfig}       
\usepackage{color}
\usepackage{fancyhdr}  
 \usepackage{natbib}  
\usepackage{rotating}
\usepackage{pdflscape}
\usepackage[T1]{fontenc}        
\usepackage{url}                
\usepackage{mathtools}                 
\usepackage{vmargin}    
\usepackage[title]{appendix}

\usepackage{lineno}
\usepackage{setspace}
\usepackage[colorlinks=true,linkcolor= blue,urlcolor=blue,citecolor=blue]{hyperref} 
\usepackage[ruled,vlined]{algorithm2e}  
\biboptions{square}   
   
\newcommand{\bo}[1]{\mathbf{#1}} 
\def\indic{\hbox{1\kern-.24em\hbox{I}}}      

\newcommand{\var}{\mathbb{V}}  
\newcommand{\esp}{\mathbb{E}}    
\newcommand{\trace}{\text{Tr}}

\newcommand{\x}{x}
\newcommand{\X}{X}

\newcommand{\Z}{Z} 
\newcommand{\R}{\mathbb{R}}
\newcommand{\proba}{\mathbb{P}} 

\newcommand{\M}{f}

\newcommand{\T}{T}    
      
\newcommand{\NN}{n}     
\newcommand{\Select}{\mathbf{s}}  
\newcommand{\rr}{\mathbf{r}} 
\newcommand{\vv}{\mathbf{v}} 
\newcommand{\oo}{\mathbf{o}} 
\newcommand{\Ori}{\boldsymbol{\pi}_{1}}
\newcommand{\norme}[1]{\left|\left| #1 \right|\right|_{L^2}}

\newcommand{\normf}[1]{\left|\left| #1 \right|\right|_F}
   
\newcommand{\normh}[1]{\left|\left| #1 \right|\right|_\mathcal{H}}

\newcommand{\D}{\Sigma}

\newcommand{\matleq}{\preceq}  

\newcommand{\vect}{\text{Vec}}   

\newtheorem{prop}{Proposition}{\bf}{\it} 
\newtheorem{defi}{Definition}{\bf}{\it}  

\newtheorem{theorem}{Theorem}{\bf}{\it}     
\newtheorem{lemma}{Lemma}{\bf}{\it}          
\newtheorem{rem}{Remark}{\bf}{\it}  
\newtheorem{corollary}{Corollary}{\bf}{\it} 

\catcode`\"=\active 
\catcode`\"=\active 
\def"{\og\ignorespaces}
\def"{{\fg}} 
 
\usepackage{etoolbox}
\makeatletter
\patchcmd{\ps@pprintTitle}{\footnotesize\itshape
       Preprint submitted to \ifx\@journal\@empty Elsevier
       \else\@journal\fi\hfill\today}{\relax}{}{}
\makeatother  
   
\begin{document}

\begin{frontmatter} 

\title{On dependency models and dependent generalized sensitivity indices}
\author[a,b]{Matieyendou Lamboni\footnote{Corresponding author: matieyendou.lamboni[at]gmail.com /@univ-guyane.fr, 18/06/2022}}          
\address[a]{University of Guyane, Department DFRST, 97346 Cayenne, French Guiana, France}
\address[b]{228-UMR Espace-Dev: University of Guyane, University of R\'eunion, IRD, University of Montpellier, France.}   
                            
\begin{abstract}  
In this paper, we derive copula-based and empirical dependency models (DMs) for simulating non-independent variables, and then propose a new way for determining the distribution of the model outputs conditional on every subset of inputs. Our approach relies on equivalent representations of the model outputs using DMs, and an algorithm is provided for selecting the representations that are necessary and sufficient for deriving  the distribution of the model outputs given each subset of inputs. In sensitivity analysis, the selected representations allow for assessing the main, total and interactions effects of each subset of inputs. We then introduce the first-order and total indices of every subset of inputs with the former index less than the latter. Consistent estimators of such indices are provided, including their asymptotic distributions. Analytical and numerical results are provided using single and multivariate response models.     
\end{abstract}       
                      
               
\end{frontmatter}        
\setpagewiselinenumbers    
\modulolinenumbers[1]   
  
\section{Introduction} 
Performing uncertainty quantification and variance-based sensitivity analysis (SA) of computer codes or mathematical models in presence of non-independent input variables remains a challenge issue when one is interested in assessing the contributions of every subset of inputs, including the interactions with other inputs. It is known that the dependence structures among the inputs may have significant impacts on the results of Sobol' indices (\cite{sobol93,sobol01}) and generalized sensitivity indices (e.g., \cite{lamboni18a, lamboni11,gamboa14,gamboa14b,xiao17,lamboni19,perrin21}). Some of such dependence structures are 
due to the imposed constraints on the model outputs and/or inputs and the correlations among inputs. \\   
 
Variance-based SA and Shapley effects rely on the distribution of the model outputs conditional on every subset of inputs, including every pair of inputs. For non-independent input variables, different approaches (see \cite{daveiga09,mara12,kucherenko12,hao13,chastaing12,kucherenko17}) provide the same first-order index of one input or one group of dependent inputs. To be able to rank the input variables using  the total indices, the works introduced in \cite{mara15,tarantola17} ensure that the first-order index is always less than the total index for each single input and for real-valued functions. In the same sense, the recent work proposed in \cite{lamboni21} provides in-depth approaches for quantifying the first-order and total effects of each single input and some subsets of inputs by making use of the dependency models of non-independent variables. Such approaches also include screening dependent inputs and dependent derivative global sensitivity measures (dDGSMs). However, these approaches do not i) allow for quantifying the effects of a pair of inputs (or more than two inputs) from the same group of dependent variables, and ii) address the issues of the model outputs or inputs subject to complex constraints.\\
    
To better highlight the motivation and the novelty of this paper in order to address the above issues, let us recall two dependency models (DMs) of the random vector $\bo{\X} :=(\X_1, \X_2, \X_3)$, that is, (\cite{skorohod76,lamboni21,lamboni22b})
$$     
(\X_2, \X_3) \stackrel{d}{=} r_1(\X_1, Z_2, Z_3);
\qquad \quad 
(\X_1, \X_3) \stackrel{d}{=} r_2(\X_2, Z_1, Z_3)  
 \, , 
$$  
where $\X_1, Z_2, Z_3$ (resp. $\X_2, Z_1, Z_3$) are independent variables. Since  DMs can be used for i) efficient simulations of dependent variables (\cite{lamboni22b}), and ii) defining dDGSMs, dependent elementary effects and dependent generalized sensitivity indices (dGSIs) (\cite{lamboni21}), the first contribution of this paper consists in deriving new dependency functions, that is, $r_j : \R^d \to \R^{d-1}$ of any $d$-dimensional random vector such as random vectors with known copula and continuous and/or discrete margins; vectors with known empirical distribution, including constrained inputs and/or outputs;  transformations of random vectors with known DMs. \\   
                  
For functions with non-independent variables (e.g., $\M(\bo{\X})$), the calculus and computations of the conditional expectations of  $\M(\bo{\X})$ given each subset of $\bo{\X}$ are relevant not only for variance based SA but also for the Shapley effects of inputs (\cite{owen14b}). For instance, the work proposed in \cite{lamboni21} combines $\M$ with the above DMs, that is,
$$
\M(\bo{\X}) \stackrel{d}{=} \M(\X_1, r_1(\X_1, Z_2, Z_3)); 
\qquad \quad     
\M(\bo{\X}) \stackrel{d}{=} \M(\X_2, r_2(\X_2, Z_1, Z_3))  \, ,  
$$
to derive the conditional expectation of $\M(\bo{\X})$ given $\X_1$ (resp. $\X_2$) and dGSIs. So far, there is no methodology for determining the conditional expectation of $\M(\bo{\X})$ given $(\X_1, \X_2)$ or $(\X_1, \X_3)$ or $(\X_2, \X_3)$ using DMs. These conditional expectations are relevant for detecting groups of the input variables that have significant impact on the model outputs. \\    

The second contribution of this paper aims at proposing a new methodology for determining the distributions of $\M(\bo{\X})$ conditional on every subset of inputs (including the conditional expectations of $\M(\bo{\X})$ given each subset of inputs) by making use of the necessary and sufficient DMs of $\bo{\X}$ or equivalent representations of $\M(\bo{\X})$. We provide an algorithm for selecting such relevant representations out of $d!$ possibilities for any $d$-dimensional random vector. We then define and study dGSIs for every subset of non-independent inputs, and introduce the asymptotic distributions of the estimators of dGSIs. \\  
       
The paper is organized as follows: Section \ref{sec:mdep} deals with generic DMs of dependent variables following any copula-based distribution. We introduce the transformation of dependency functions, which avoids searching the DMs of some random vectors by making use of known dependency functions. We also derive empirical dependency functions when the analytical distributions of input variables are not available. In Section \ref{sec;repfun}, we introduce equivalent representations of the model outputs that share the same distribution and are relevant for determining the conditional expectations of the model outputs given each subset of inputs. Section \ref{sec:msa} deals with the extension of  dependent multivariate SA (dMSA) (\cite{lamboni21}) and the estimators of the new dGSIs. We also provide the asymptotic distributions of dGSIs. Analytical and numerical results are provided in Section \ref{sec:test}, and we  conclude this work in Section~\ref{sec:con}.
						      		  			            		              
\section*{General notation}                      
For integer $d >0$, we use $\bo{\X} :=(X_1, \, \ldots,\, X_d)$ for a random vector.  For $u \subseteq \{1,\,\ldots,\, d\}$, we use $\bo{\X}_u :=(\X_j, \forall\, j \in u)$, $\bo{\X}_{\sim u} :=(\X_j, \forall\, j\in \{1,\,\ldots,\, d\}\setminus u)$ and $|u|$ for its cardinality (i.e., the number of elements in $u$). Thus, we have the partition $\bo{\X} =(\bo{\X}_{u},\, \bo{\X}_{\sim u})$. We also use $\bo{Z} \stackrel{d}{=} \bo{\X}$ to say that $\bo{Z}$ and $\bo{\X}$ have the same distribution (CDF). For $\bo{a} \in\R^\NN$, we use $\norme{\bo{a}}$ for the Euclidean norm. For a matrix $\Sigma \in \R^{\NN\times \NN}$, we use $\trace(\Sigma)$ for the trace of $\Sigma$, and $\normf{\Sigma} :=\sqrt{\trace\left(\Sigma \Sigma^\T \right)}$ for the Frobenius norm of $\Sigma$. 
We use $\esp[\cdot]$ for the expectation and $\var[\cdot]$ for the variance-covariance.
                        
\section{Dependency functions of non-independent random variables} \label{sec:mdep}  
This section provides generic DMs of random vectors following distributions that were not considered in \cite{lamboni21,lamboni22b}, including empirical distributions and the DMs of some transformations of random vectors with known dependency functions. \\     
    
Namely, the inputs $\bo{\X}$ have $F$ as the joint CDF, $C$ as copula, and  $F_j$ as the marginal CDF of $\X_j$ with $j=1,\ldots, d$.  For all $\bo{x} \in \R^d$, the joint CDF at $\bo{x}$ is given by $F(\bo{x})=C(F_1(x_1),\, \ldots,\, F_d(x_d))$. We use $F_j^{\leftarrow}$ for the generalized inverse of $F_j$, and  $F_{j|k}$ for the distribution of $X_j$ conditional on $X_k$ for all $j, \, k \in \{1, \ldots, d\}$. For a discrete variable $\X_i$ having $F_i$ as the marginal CDF and $\lambda_i \in [0, 1]$, the distributional transform of $\X_i$ is given by (\cite{ferguson67,ruschendorf81,ruschendorf05,ruschendorf09}) 
$$
\tau_{F_i}(x_i, \lambda_i) := \proba(\X_i < x_i) + \lambda_i \proba(\X_i =x_i)
$$ 
and it ensures that for all $U_i  \sim \mathcal{U}(0,\, 1)$ 
\begin{equation} \label{eq:dt1} 
V_i :=\tau_{F_i}(\X_i, U_i) \sim  \mathcal{U}(0,\, 1) \, ,    
\qquad 
	X_i \stackrel{d}{=} \tau_{F_i}^{\leftarrow}(V_i) \, ,
	\quad
X_i \stackrel{d}{=} F_i^{\leftarrow} (V_i) \;  
	\; \mbox{a.s.} \, .
\end{equation} 
It is to be noted that the above distributional transform comes down to the Rosenblatt transform (\cite{rosenblatt52}) for a continuous CDF (i.e., $F_j$), that is, $\tau_{F_j}(\X_j, U_j) =F_j(\X_j) \sim  \mathcal{U}(0,\, 1)$. Also, remark that $X_i \stackrel{d}{=} F_i^{\leftarrow} (V_i) $ is  equivalent to the inverse of the Rosenblatt transform, which is part of the multivariate conditional quantile transform (MCQT) (\cite{obrien75,arjas78,ruschendorf81}). The MCQT implies a regression representation of $\bo{\X}$ (\cite{skorohod76,ruschendorf93,ruschendorf05,ruschendorf09}), which also implies a dependency model of $\bo{\X}$ (\cite{skorohod76,lamboni21,lamboni22b}).\\
   
To recall the theoretical expression of a DM in Equation (\ref{eq:depmgen}), we use $(w_1, \ldots, w_{d-1})$ for an arbitrary permutation of $\{1, \ldots, d\}\setminus\{j\}$, $\bo{\X}_{\sim j} := (\X_{w_1}, \ldots, \X_{w_{d-1}})$ and $\bo{Z} := (Z_{w_1}, \ldots, Z_{w_{d-1}}) \sim \mathcal{U}(0,\, 1)^{d-1}$. A dependency function $r_j := (r_{w_{1}}, \ldots, r_{w_{d-1}}) : \R^d \to \R^{d-1}$ is a measurable function that satisfies (\cite{skorohod76,lamboni21,lamboni22b}): 
\begin{eqnarray}  \label{eq:depmgen}   
\bo{\X}_{\sim j} \stackrel{d}{=}  r_j \left(\X_j, \bo{Z} \right)
  =  \left[ \begin{array}{l} 
 r_{w_1} \left(\X_j,\, Z_{w_1} \right) \\ 
  r_{w_2} \left(\X_j,\, Z_{w_1}, \, Z_{w_2} \right) \\ 
  \vdots \\
 r_{w_{d-1}} \left(\X_j,\, \bo{Z} \right) 
\end{array}  \right];
\; \quad \mbox{and} \quad   (\X_j, \bo{\X}_{\sim j}) \stackrel{d}{=} \left(\X_j, r_j \left(\X_j, \bo{Z} \right) \right)\, ,                      
\end{eqnarray}         
where $\bo{Z}$  and $\X_j$ are independent. It is worth noting that $r_j$ is not unique. Indeed, one may replace in Equation (\ref{eq:depmgen}) $Z_{w_i}$ with $F_{R_{w_i}}(R_{w_i})$ for any continuous variable $R_{w_i} \sim F_{R_{w_i}}$ or with $\tau_{F_{W_{w_i}}}(W_{w_i}, U_{w_i})$ for any discrete variable $W_{w_i} \sim F_{W_{w_i}}$.  Immediate and interesting properties of the DM provided in Equation (\ref{eq:depmgen}) are gven below.  
          
\begin{prop} \label{prop:gdm1}
For an integer $p \leq d-1$, let $u =(j,\, w_1, \ldots, w_p)$ be a vector. Then, there exists $r_{u} : \R^d \to \R^{d-|u|}$ such that 
\begin{equation} \label{eq:extmodel}  
\bo{\X}_{u} \stackrel{d}{=} \left(\X_j, \, r_{w_1} \left(\X_j,\, Z_{w_1} \right), \ldots, r_{w_p} \left(\X_j, \, Z_{w_1}, \ldots, Z_{w_{p}} \right) \right) \, ;    
\end{equation}  
and  
\begin{equation}   \label{eq:extmodelper}
\bo{\X}_{\sim u} \stackrel{d}{=}  
r_{u}\left(\bo{\X}_{u}, Z_{w_{p+1}}, \ldots, Z_{w_{d-1}} \right)   \, .       
\end{equation}
\end{prop}        
The proof is straightforward using the regression representation of random vectors (\cite{ruschendorf09}) and DMs.  
 
\subsection{Distribution-based and copula-based  dependency functions} \label{sec:expdepfun} 
Using the joint and marginal CDFs of $\bo{\X}$, the distribution-based dependency model of $\bo{\X}$ is given by (\cite{skorohod76,lamboni21})  
\begin{eqnarray}        \label{eq:depdf}
\begin{array}{ccl}               
\X_{w_1}   & = & r_{w_1} \left(\X_j,\, Z_{w_1} \right) = F_{w_1 |j}^{\leftarrow} \left(Z_{w_1} \,|\, \X_j \right) \\              
\X_{w_{d-1}} &=&  F_{w_{d-1} |j, w_1,\ldots, w_{d-2}}^{\leftarrow} \left(Z_{w_{d-1}} \,|\,\X_j, \, r_{w_1}(\X_j, \,Z_{w_1}), \ldots, r_{w_{d-2}}(\X_j, \,Z_{w_1},\ldots, Z_{w_{d-2}}) \right) \\ 
			\end{array}    \, ,   
\end{eqnarray}   
       
Likewise, the copula-based dependency function is of great interest to deal with all the distributions having the same copula $C$ (\cite{nelsen06,ruschendorf05}). When all the marginal CDFs are continuous, the copula-based expression of a DM is derived as follows: 
\begin{eqnarray}     \label{eq:depcopcont}   
\begin{array}{ccl}                 
\X_{w_1}   & = & r_{w_1} \left(\X_j, Z_{w_1} \right) =  F_{w_{1}}^{\leftarrow} \left( C_{w_1 | j}^{\leftarrow}\left(Z_{w_1} \, | F_j(\X_j) \right) \right) \\  
\X_{w_2}   & = & r_{w_2} \left(\X_j, Z_{w_1}, Z_{w_2} \right) =  F_{w_{2}}^{\leftarrow} \left( C_{w_2 | j, w_1}^{\leftarrow}\left(Z_{w_2} \, | F_j(\X_j), F_{w_1}\left(r_{w_1} \left(\X_j, Z_{w_1} \right) \right) \right) \right)\\        
\vdots & & \\ 
\X_{w_{d-1}} &=& F_{w_{d-1}}^{\leftarrow} \left( C_{w_{d-1} |j, w_1,\ldots, w_{d-2}}^{\leftarrow}\left(Z_{w_{d-1}} \, |\, F_j(\X_j), F_{w_1}\left(r_{w_1} \left(\X_j, Z_{w_1} \right) \right), \ldots \right) \right)   \\  
			\end{array}    \, ,               
\end{eqnarray}          
bearing in mind the general sampling algorithm based on copulas (\cite{nelsen06,mcneil15}). To cope with discrete variables, we consider independent variables, that is,  $\bo{U} \sim \mathcal{U}(0,\, 1)^d$, $\bo{Z}$ and $\X_j$, and we use $V_{j} := (\X_j, U_j)$; $V_{w_k} := (Z_{w_k}, U_{w_k})$ for all $k \in \{1,\ldots, d\} \setminus \{j\}$. The general copula-based DM is then derived as follows: 
   
\begin{eqnarray}     \label{eq:depcop}   
\begin{array}{ccl}                 
\X_{w_1}   & = & r_{w_1} \left(V_j, Z_{w_1} \right) =  F_{w_{1}}^{\leftarrow} \left( C_{w_1 | j}^{\leftarrow}\left(Z_{w_1} \, | \tau_{F_j}(V_j) \right) \right) \\   
\X_{w_2}   & = & r_{w_2} \left(V_j, V_{w_1}, Z_{w_2} \right) =  F_{w_{2}}^{\leftarrow} \left( C_{w_2 | j, w_1}^{\leftarrow}\left(Z_{w_2} \, | \tau_{F_j}(V_j), \tau_{F_{w_{1}}} \left( r_{w_{1}}(V_j, Z_{w_1}), U_{w_1}\right) \right) \right) \\  
 \vdots & & \\  
\X_{w_{d-1}} &=& F_{w_{d-1}}^{\leftarrow} \left( C_{w_{d-1} |j, w_1,\ldots, w_{d-2}}^{\leftarrow}\left(Z_{w_{d-1}} \, |\, \tau_{F_j}(V_j), \tau_{F_{w_{1}}} \left( r_{w_{1}}(V_j, Z_{w_1}), U_{w_1}\right), \ldots \right) \right)   \\  
			\end{array}    \, .  
\end{eqnarray}             
 
For the Gaussian copula and continuous marginal CDFs, the dependency function in (\ref{eq:depcopcont}) has a simple expression, and it has been provided in \cite{lamboni21}. Lemma \ref{lem:gaus} extends such function to cope with discrete and continuous variables. To that end, consider the Gaussian copula  given by $C^{Gauss}(U_1, \ldots, U_d,\, \mathcal{R})$ with $\mathcal{R}$ the correlation matrix. We use $\mathcal{L}$ for the Cholesky factor of $\mathcal{R}$ (i.e., $\mathcal{R} =\mathcal{L}\mathcal{L}^\T$), $\Phi$ for the CDF of the standard Gaussian variable and $\mathcal{I}$ for the identity matrix.
  
\begin{lemma} \label{lem:gaus}    
Let $\X_j$, $\bo{Z} \sim \mathcal{N}_{d-1} \left(\bo{0}, \mathcal{I} \right)$ and $U_j \sim \mathcal{U}(0, 1)$ be independent variables.\\     

If $(\X_{j},\, \bo{\X}_{\sim j})$ has $C^{Gauss}\left(U_j,\, U_{w_1},\ldots, U_{w_{d-1}}, \mathcal{R}\right)$ as copula, then
\begin{equation} \label{eq:depgauco}    
   \bo{\X}_{\sim j}  \stackrel{d}{=}  r_j\left(\X_{j},\bo{\Z}, U_j \right) =
\left[\begin{array}{ccl}    
\X_{w_1}   & = & F_{{w_1}}^{\leftarrow}\left[\Phi \left(Y_{w_1} \right)\right] \\              
  &\vdots &  \\           
\X_{w_{d-1}} &=&  F_{w_{d-1}}^{\leftarrow}\left[\Phi \left(Y_{w_{d-1}} \right)\right] \\
							\end{array}                       
								\right]          \, ,                  
\end{equation}   
with  $$    
\left[\begin{array}{c}
Y_j \\
\bo{Y}_{\sim j}
\end{array} \right]  
:=\mathcal{L}  \left[
\begin{array}{c}
\Phi^{-1}\left(\tau_{F_j}(X_j, U_j)\right) \\  
 \bo{\Z}  \\ 
\end{array} 
\right]   \, .  
$$ 
\end{lemma}
 \begin{preuve}      
See Appendix \ref{app:lem:gaus}.
\begin{flushright}         
$\Box$        
\end{flushright} 
\end{preuve}     
 
To provide the dependency function for the Student copulas in Lemma~\ref{lem:stucop}, we use $t(\nu,\, 0,\, 1)$ for the standard $t$-distribution with $\nu$ degrees of freedom and $T_\nu$ for its CDF; $C^{St}\left(U_j,\, U_{w_1},\ldots, U_{w_{d-1}}, \, \nu,\,\mathcal{R}\right)$ for the Student copula. 
\begin{lemma} \label{lem:stucop}
Let $\X_j$, $\left\{Z_{w_i} \sim t(\nu + i,\, 0,\, 1)\right\}_{i=1}^{d-1}$ and $U_j$ be independent variables.\\ 
 
 If $(\X_{j},\, \bo{\X}_{\sim j})$  has $C^{St}\left(U_j,\, U_{w_1},\ldots, U_{w_{d-1}}, \, \nu,\,\mathcal{R}\right)$ as copula, then  
\begin{equation} \label{eq:depgauco}   
   \bo{\X}_{\sim j}  \stackrel{d}{=}  r_j\left(\X_{j},\bo{\Z}, U_j \right) =
\left[\begin{array}{ccl}    
\X_{w_1}   & = & F_{{w_1}}^{\leftarrow}\left[T_\nu \left(Y_{w_1} \right)\right] \\            
  &\vdots &  \\           
\X_{w_{d-1}} &=&  F_{w_{d-1}}^{\leftarrow}\left[T_\nu \left(Y_{w_{d-1}} \right)\right] \\
							\end{array}                       
								\right]          \, ,                  
\end{equation}     
with 
$$ \displaystyle   
\left[\begin{array}{c}
Y_j \\ 
\bo{Y}_{\sim j}
\end{array} \right] 
  \stackrel{}{:=}    
 \mathcal{L} \; 
         \left[\begin{array}{c} 
				T_\nu^{-1} \left(\tau_{F_j}(X_j, U_j)\right)  \\ 
	\sqrt{\frac{\nu + \left(T_\nu^{-1} \left(\tau_{F_j}(X_j, U_j)\right) \right)^2}{\nu +1}} Z_{w_1} \\       
		\vdots \\
	\sqrt{\frac{\left(\nu + \left(T_\nu^{-1} \left(\tau_{F_j}(X_j, U_j)\right)\right)^2 \right) \prod_{k=1}^{d-2}\left(\nu + k+ Z_{w_k}^2 \right)}{\prod_{k=1}^{d-1}(\nu +k)}} Z_{w_{d-1}} \\
							\end{array} \right]  \, .  
$$         
\end{lemma}  
 \begin{preuve}      
See Appendix \ref{app:lem:stucop}. 
\begin{flushright}         
$\Box$        
\end{flushright}  
\end{preuve}        
 
It comes out that we have to include additional and independent  uniformly distributed variables in the copula-based dependency functions in presence of discrete variables. It is worth noting that we must replace $\tau_{F_j}(X_j, U_j)$ with $F_j(X_j)$ in Lemmas~\ref{lem:gaus}-\ref{lem:stucop} for every continuous variable $\X_j \sim F_j$. 
           
\subsection{Transformation of dependency functions} \label{sec:expdepftrans}
Many transformations such as monotonic transformations are widely used in probability and statistics. In this section, we investigate the derivation of the dependency functions of dependent variables defined through such transformations. \\

For an invertible transformation $T_j$, $T_j^{\leftarrow}$ denotes its inverse.
 We use $R$ for the Rademacher random variable, that is, the probability $\proba(R=1)=\proba(R=-1)=1/2$;  $\bo{R}:=(R_i \sim R, i=1, \ldots, \NN)$ for independent variables and $diag(\bo{R}) \in \R^{\NN \times \NN}$ for a diagonal matrix.  
 
\begin{prop} \label{prop:trans}         
Let $r_j$ be a dependency function of $\bo{\X}$, $T_{j} :\R \to \R$ be an invertible transformation and $T_{\sim j} : \R^{d-1} \to \R^{\NN}$ be a measurable function. \\ 
    
$\quad (i)$ If $(Y_{j},\, \bo{Y}_{\sim j}) \stackrel{d}{=} \left(T_{j}(\X_{j}), \, T_{\sim j}(\bo{\X}_{\sim j}) \right)$, then a DM of $\bo{Y}$ is given by 
\begin{equation}  \label{eq:trans} 
\displaystyle      
\bo{Y}_{\sim j} = T_{\sim j}\left(r_{j}\left(T_{j}^{\leftarrow}(Y_{j}), \, \bo{Z} \right) \right) \, .        
\end{equation}   
	   
$\quad (ii)$ If $\left(|Y_{j}|, \, |\bo{Y}_{\sim j}|\right) \stackrel{d}{=} \left(T_{j}(\X_{j}), \, T_{\sim j}(\bo{\X}_{\sim j}) \right)$ and every $Y_{w_i}$ has a symmetric distribution about $0$, then 
\begin{equation} \label{eq:abstrans} 
\displaystyle     
\bo{Y}_{\sim j} \stackrel{\text{}}{=} diag(\bo{R}) \; T_{\sim j} \left( r_{j}\left(T_{j}^{\leftarrow}\left(|Y_{j}|\right), \, \bo{Z}  \right)\right) \, .
\end{equation}  
\end{prop}         
\begin{preuve}            
See Appendix \ref{app:prop:trans}.  
\begin{flushright}
$\Box$ 
\end{flushright}
\end{preuve}       
 
\begin{corollary} \label{coro:trans}        
Let $T_k,\, k =1, \ldots, d$ be measurable functions and $T_{j}$ be invertible. \\
         
$\quad (i)$ If $Y_k  =T_k(\X_k)$ with $k =1, \ldots, d$, then $\forall\, i \in\{1,\ldots, d-1\}$, 
\begin{equation}  \label{eq:transindiv}   
\displaystyle      
Y_{w_i} = T_{w_i}\left(r_{w_i}\left(T_{j}^{\leftarrow}(Y_{j}), \, \Z_{w_1}, \ldots, Z_{w_i} \right) \right) \, .  
\end{equation} 

 $\quad (ii)$ If $|Y_k| = T_k(\X_k)$ and $Y_{w_i}$ has a symmetric distribution about $0$, then    
 \begin{equation}  
 \displaystyle   
 Y_{w_i} \stackrel{}{=} R_{w_i} \times T_{w_i} \left( r_{w_i}\left(T_j^{\leftarrow}\left(|Y_{j}|\right), \, \Z_{w_1}, \ldots, Z_{w_i} \right)\right) \, . 
 \end{equation}      
\end{corollary}  
  
When the transformation $T_j$ is strictly increasing on the support of $\X_j$ for all $j \in \{1, \ldots, d\}$, Equation (\ref{eq:transindiv}) provides a dependency function of  $\left(T_1(\X_1), \ldots, T_d(\X_d) \right)$. For continuous variables, it is known that $\left(T_1(\X_1), \ldots, T_d(\X_d) \right)$ and $\bo{\X}$ have the same copula but different margins. Therefore, Corollary \ref{coro:trans} gives a new way for deriving the DM of the random vector $\left(T_1(\X_1), \ldots, T_d(\X_d) \right)$ without explicitly using its marginal CDFs. Thus, Proposition \ref{prop:trans} and Corollary \ref{coro:trans} extend the multivariate conditional method, as DMs can be used for sampling random values as well.  

\subsection{Empirical dependency functions}
This section deals with the derivation of DMs for unknown distributions of inputs such as distributions obtained by imposing constraints on the initial inputs or outputs. Formally, given a function $c : \R^d \to \R^\NN$ and a domain of interest $D$, we are interested in deriving a dependency function of a random vector defined by  
\begin{equation} \label{eq:dfcons} 
\bo{\X}^c \stackrel{d}{=} \left\{\bo{\X} \sim F :\, c(\bo{\X}) \in D \right\} \, .  
\end{equation}  
   
While we are able to derive the analytical distribution of $\bo{\X}^c$ and its dependency function for some distributions and constraints (see \cite{lamboni21,lamboni22b}), we have to estimate such dependency functions in general. Using Equation (\ref{eq:dfcons}), we can generate a sample of $\bo{\X}^c$, that is,  $\bo{\X}^c_1, \ldots, \bo{\X}^c_m$ and a pseudo-sample from the copula $C$ of $\bo{\X}^c$, that is, $\hat{F}_{1}(\X_{i, 1}^c), \ldots, \hat{F}_{d}(\X_{i,d}^c)$ with $i=1, \ldots, m$ for continuous variables, where $\hat{F}_j$ is an estimator of $F_j$. In general, we are going to use  $\hat{\tau_{F_1}}(\X_{i, 1}^c,  U_1), \ldots, \hat{\tau_{F_d}}(\X_{i,d}^c, U_d)$. \\   
    
With such samples, we consider two main ways to derive the empirical dependency functions. Firstly, we fit a distribution to that observations, and then derive dependency functions using results from Section \ref{sec:mdep}. To that end, there are numerous papers about fitting a distribution to data. For instance, direct methods for estimating densities and distributions can be found in \cite{rosenblatt56,parzen62,epanechnikov69,silverman86}, and the copula-based methods for modeling distributions are provided in (\cite{clayton78,joe97,mcneil15,ruschendorf05,smith10,durante15}). \\   

Secondly, we derive empirical dependency functions using the estimators of the conditional quantile functions (\cite{koenker78,truong89,hendricks92,koenker94,koenker01,bassett82,koenker05,takeuchi06}). Formally, consider the loss function of Koenker and Bassett (\cite{koenker78}) given by $L(x, u) = x\left(u-\indic_{\{x < 0\}} \right)$ with $u \in [0,\, 1]$ and $\indic_{\{x < 0\}}$ the indicator function. A dependency function can be written as follows:
\begin{equation} \label{eq:mestdepf}
r_j(\X_j^c,\, Z_{w_1}) := \arg\min_{f \in \mathcal{F}} \,\esp\left[ L(\X_{w_1}^c-f(\X_j^c),  Z_{w_1})\, | \, \X_j^c,\, Z_{w_1} \right] \, ,    \nonumber
\end{equation}
where $\mathcal{F}$ is a class of smooth functions, and $Z_{w_1} \sim \mathcal{U}(0,\, 1)$ is independent of $(\X_j^c,\, \X_{w_1}^c)$. 
 Using the sample of $(\X_j^c,\, \X_{w_1}^c)$, the M-estimator of a dependency function is given by (\cite{takeuchi06}, Lemma 3)
\begin{equation} 
\widehat{r_j}(\X_j^c,\, Z_{w_1}) := \arg\min_{f \in \mathcal{H}} \sum_{i=1}^m L(\X_{i, w_1}^c-f(\X_{i, j}^c),  Z_{w_1}) + \frac{\lambda}{2} \normh{f(\X_{j}^c) -b}^2 \, , 
\end{equation}
where $\lambda \in \R_+$ is a bandwidth, $\normh{\cdot}$ is a given norm, and $b\in \R$.
        
\begin{rem}  
When the constraints involve a model that is time demanding to obtain one run,  we may replace such model with its emulator. The Gaussian process (\cite{currin91,haylock96,kennedy01,jeremy04}) is an interesting candidate, as it deals with dependent variables. 
\end{rem}  
     
\section{Representation of functions with dependent variables}  \label{sec;repfun}
This section formalizes different representations of a function with non-independent variables, which are relevant for determining the conditional expectations of the model outputs given any subset of the inputs.  
Formally, consider a vector-valued function  $\M : \R^d \to \R^\NN$ that includes $d$ random inputs $\bo{\X}$ and provides $\M(\bo{\X})$ as outputs. The input variables $\bo{\X}$ can always be organized as follows (see also \cite{lamboni21}): \\  

(A1): the random vector $\bo{\X} = \left(\X_1, \ldots, \X_d \right)$ is consisted of $K\geq 1$ independent random vectors, that is, $\bo{\X} := \left(\bo{\X}_{\boldsymbol{\pi}_1}, \ldots, \bo{\X}_{\boldsymbol{\pi}_K}\right)$
where   
$ 
\bo{\X}_{\boldsymbol{\pi}_k} = \left(\X_{1,k}, \ldots, \X_{d_k,k} \right)$
with $k=1, \ldots, K$, and $\bo{\X}_{\boldsymbol{\pi}_{k_1}}$, $\bo{\X}_{\boldsymbol{\pi}_{k_2}}$ are independent for all $k_1 \neq k_2$. Thus, $\boldsymbol{\pi}_{k}$ is the set of all the subscripts of the variables that are in $\bo{\X}_{\boldsymbol{\pi}_k}$, that is, $\bo{\X}_{\boldsymbol{\pi}_k} =: \left( \X_i,\, \forall\, i \in \boldsymbol{\pi}_k \right)$.    
Without loss of generality, we use $\bo{\X}_{\boldsymbol{\pi}_{1}}$ for a vector of $d_1\geq 0$ independent initial variable(s) and  $\bo{\X}_{\boldsymbol{\pi}_{k}}$ for a vector of $d_k\geq 2$ dependent variables with $k \geq 2$.\\ 
    
In what follows, $\bo{w}_k := (w_{1,k},\ldots, w_{d_k,k})$ denotes an arbitrary permutation of the elements of $\boldsymbol{\pi}_{k}$ and $\bo{w}_{\sim 1, k} := (w_{2,k},\ldots, w_{d_k,k})$. Thus, $w_{i,k} \in \boldsymbol{\pi}_{k},\; \forall \, i \in  \{1, \ldots, d_k\}$. If we use  $\Select := \{w_{1,2}, \ldots, w_{1,K}\}$, then $\bo{\X}_{\Select} := \left(\X_{w_1,2}, \ldots, \X_{w_1,K} \right)$ contains $K-1$ independent variables and $\bo{\X}_{\sim \Select} :=\left(\bo{\X}_{\bo{w}_{\sim 1, 2}}, \ldots, \bo{\X}_{\bo{w}_{\sim 1, K}} \right)$ with $\bo{\X}_{\bo{w}_{\sim 1, k}} := \left(\X_{w_2,k}, \ldots, \X_{w_{d_k},k} \right)$, $k=2, \ldots, K$.
For an integer $p_k$ with $p_k \leq d_k$ and $k=2, \ldots, K$, let us define 
\begin{equation} \label{eq:setvuk}
\rr_k := (w_{2,k},\ldots, w_{p_k,k}); \qquad \vv_k := (w_{1,k},\, \rr_k): 
\qquad 
 \boldsymbol{\varpi}_k :=(w_{p_k+1,k}, \ldots, w_{d_k,k}) \, ,  
\end{equation}      
By definition, $\bo{w}_{\sim 1, k} =(\rr_k, \, \boldsymbol{\varpi}_k)$; $\rr_k$ and $\vv_k$ vanish when $p_k=0$, and $\vv_k=w_{1,k}$ when $p_k=1$.        
A Dependency model of $\bo{\X}$ is given by (see Section \ref{sec:mdep} and \cite{lamboni21,lamboni22b})
\begin{equation} \label{eq:mdepk}       
\bo{\X}_{\sim \Select} = \left(
\bo{\X}_{\bo{w}_{\sim 1, 2}} \stackrel{d}{=} r_{w_{1, 2}}\left(\X_{w_{1, 2}},\bo{Z}_{\bo{w}_{\sim 1, 2}} \right),
  \ldots ,
 \bo{\X}_{\bo{w}_{\sim 1, K}}  \stackrel{d}{=} r_{w_{1, K}}\left(\X_{w_{1, K}},\bo{Z}_{\bo{w}_{\sim 1, K}} \right)
\right) \,   ;  
\end{equation} 
where $\bo{Z}_{\bo{w}_{\sim 1, k}} :=\left(Z_{w_{2, k}}, \ldots, Z_{w_{d_k}, k} \right)$ contains independent variables; $\X_{w_{1,k}}$  is independent of $\bo{Z}_{\bo{w}_{\sim 1, k}}, \, k=2, \ldots, K$. If we use $\bo{Z}_{\bo{w}_{\sim 1}} :=(\bo{Z}_{\bo{w}_{\sim 1, 2}},\ldots, \bo{Z}_{\bo{w}_{\sim 1, K}})$, then $(\bo{\X}_{\Ori}, \,\bo{\X}_{\Select},\, \bo{Z}_{\bo{w}_{\sim 1}})$ contains independent variables under (A1). Equation (\ref{eq:mdepk}) can be written in a concise way using the function $r_{\Select} : \R^{d-d_1} \to \R^{d-d_1-K+1}$, that is,
\begin{equation}    
\bo{\X}_{\sim \Select} \stackrel{d}{=} r_{\Select}\left(\bo{\X}_{\Select}, \bo{Z}_{\bo{w}_{\sim 1}} \right) =\left(r_{w_{1, 2}}\left(\X_{w_{1, 2}},\bo{Z}_{\bo{w}_{\sim 1, 2}} \right), 
  \ldots,   
	r_{w_{1, K}}\left(\X_{w_{1, K}},\bo{Z}_{\bo{w}_{\sim 1, K}} \right)
\right) \, , \nonumber 
\end{equation}  
and the properties of DMs allow for writing 
$
\bo{\X} \stackrel{\text{d}}{=}   \left(\bo{\X}_{\Ori},\, \bo{\X}_{\Select},\, r_{\Select}\left(\bo{\X}_{\Select}, \bo{Z}_{\bo{w}_{\sim 1}} \right) \right)$
and  
$
\M(\bo{\X} ) \stackrel{\text{d}}{=} 
  \M\left(\bo{\X}_{\Ori}, \bo{\X}_{\Select},\, r_{\Select}\left(\bo{\X}_{\Select}, \bo{Z}_{\bo{w}_{\sim 1}} \right) \right)
$ (\cite{lamboni21}).  
Consider the function   
\begin{equation} \label{eq:indeg} 
g\left(\bo{\X}_{\Ori}, \bo{\X}_{\Select}, \bo{Z}_{\bo{w}_{\sim 1}} \right) := \M\left(\bo{\X}_{\Ori}, \bo{\X}_{\Select},\, r_{\Select}\left(\bo{\X}_{\Select}, \bo{Z}_{\bo{w}_{\sim 1}} \right) \right) \, , 
\end{equation}                
which includes only independent variables, that is, $\left(\bo{\X}_{\Ori}, \bo{\X}_{\Select}, \bo{Z}_{\bo{w}_{\sim 1}}\right)$. Lemma \ref{lem:equivdist} provides useful properties of $g$ that are linked to $\M$. It is common to use $\M\left(\bo{\X} \right) |\, \bo{\X}_{u}$ for a random vector that represents the outputs $\M\left(\bo{\X} \right)$ conditional on $\bo{\X}_{u}$.
 
\begin{lemma} \label{lem:equivdist}     
Let $\vv_1 \subseteq \Ori$, $\{k_1, \ldots, k_{m} \} \subseteq  \{2, \ldots, K \}$, $\{w_{1,k_1}, \ldots, w_{1,k_{m}} \} \subseteq  \Select$. Then,
\begin{equation}        \label{eq:equivg}        
 \M\left(\bo{\X} \right) |\, \bo{\X}_{\vv_1}, \bo{\X}_{\vv_{k_1}}, \ldots, \bo{\X}_{\vv_{k_m}}  \stackrel{d}{=}  g\left(\bo{\X}_{\Ori}, \bo{\X}_{\Select}, \bo{Z}_{\bo{w}_{\sim 1}} \right) | \, \bo{\X}_{\vv_1}, \X_{w_{1,k_1}}, \bo{Z}_{\rr_{k_1}}, \ldots, \X_{w_{1,k_m}}, \bo{Z}_{\rr_{k_m}} \, .  \nonumber   
\end{equation}                                        
\end{lemma}      
\begin{preuve} 
See Appendix \ref{app:lem:equivdist}. 
\begin{flushright}
$\Box$          
\end{flushright} 
\end{preuve}

It comes out from Lemma \ref{lem:equivdist} that the distribution of $\M(\bo{\X})$ conditional on $\bo{\X}_u := \left(\bo{\X}_{\vv_1},\, \bo{\X}_{\vv_{k_1}}, \ldots, \bo{\X}_{\vv_{k_m}}\right)$ is equal to the distribution of $g\left(\bo{\X}_{\Ori}, \bo{\X}_{\Select}, \bo{Z}_{\bo{w}_{\sim 1}} \right)$ conditional on $\left(\bo{\X}_{\vv_1},\, \X_{w_{1,k_1}}, \bo{Z}_{\rr_{k_1}}, \ldots, \X_{w_{1,k_m}}, \bo{Z}_{\rr_{k_m}} \right)$. Thus, we are able to assess the effects of $\bo{\X}_u$ on $\M(\bo{\X})$ using $g\left(\bo{\X}_{\Ori}, \bo{\X}_{\Select}, \bo{Z}_{\bo{w}_{\sim 1}} \right)$ and $\left(\bo{\X}_{\vv_1}, \X_{w_{1,k_1}},  \bo{Z}_{\rr_{k_1}}, \ldots, \X_{w_{1,k_m}}, \bo{Z}_{\rr_{k_m}} \right)$. This leads to the following definition. 
 
\begin{defi} \label{def:equirep}
Consider $u \subseteq \{1, \ldots, d\}$ and $g$ given by Equation (\ref{eq:indeg}). \\
  
 A function $g$ is said to be an equivalent representation of $\M$ regarding the input(s) $\bo{\X}_{u}$ if the distribution of $\M(\bo{\X}) |\bo{\X}_u$ can be determined using $g$ and some of its inputs.
\end{defi}  
 
\noindent 
Different and equivalent representations of $\M$ (i.e., $g$) are necessary for assessing the effects of $\bo{\X}_u$ on the outputs for all $u \subseteq \{1, \ldots, d\}$. For instance, the representation in Lemma \ref{lem:equivdist} can be used to assess the effects of the subset of inputs given by 
\begin{eqnarray} 
\left\{\left(\bo{\X}_{\vv_1}, \bo{\X}_{\vv_{k_1}}, \ldots, \bo{\X}_{\vv_{k_m}} \right) \, : 
\begin{array}{l}   
\forall\, \vv_1 \subseteq \Ori,
\forall\, \{k_1, \ldots, k_m\} \subseteq \{2, \ldots, K\}, \\  
\forall\, p_{k_i} \in \{0, \ldots, d_{k_i}\},\, i=1, \ldots, m   \\ 
\end{array}  
 \right\}   \nonumber \, .    
\end{eqnarray} 
Modifying $\Select$ and $\rr_k$ with $k=2,\ldots, K$ leads to another representation of $\M$, which allows for assessing the effects of other inputs such as $\X_{w_\imath,k}$ with $w_{\imath,k} \notin \Select$. A permutation of the elements of $\boldsymbol{\pi}_{k}$ gives such modifications, and  we have $(d_k-1)! \left( \prod_{\substack{i=2\\ i\neq k}}^K d_i! \right)$ equivalent representations of $\M$ that can be used for determining the effects of $\X_{w_\imath,k}$. 
     
\begin{defi}
Let $u \subseteq \{1, \ldots, d\}$ and $g_1\neq g_2$ be two equivalent representations of $\M$. \\ 

The representations $g_1,\, g_2$ are said to be replicated representations of $\M$ regarding $\bo{\X}_u$ if  $g_1,\, g_2$ allow for determining the distribution of $\M(\bo{\X}) |\bo{\X}_u$.  
\end{defi}     
        
Obviously, replicated representations should be avoided when we are only interested in the distribution of $\M(\bo{\X})$ given $\bo{\X}_u$. However, note that such representations are sometime necessary to determine the distributions of the outputs conditional on other inputs. Therefore, replicated representations regarding an input or group of inputs should be avoided as much as possible. Thus, we are going to use an algorithm to select the equivalent representations of $\M$ that are necessary and sufficient for assessing the effects of all the subsets of inputs. \\ 

Keeping in mind Lemma \ref{lem:equivdist}, we need to select the permutations of $\boldsymbol{\pi}_{k}$ of the form $(w_{1,k}, \ldots, w_{d_k,k})$ in such a way that  $\forall \, u_k \subseteq \boldsymbol{\pi}_{k}$, there exists $p_k \leq d_k$ such that $u_k =\{w_{1,k}, \ldots, w_{p_k,k}\}$  with $k=2,\ldots, K$ (see Equation (\ref{eq:setvuk})). Algorithm \ref{algo:sets} aims at finding such permutations and provides the minimum number of the permutations needed. For instance, for a given permutation $(w_{1,k}, \ldots, w_{d_k,k})$, we obtain the following subsets:  
$$ 
u_k=\{w_{1,k}\}, \{w_{1,k}, w_{2,k}\}, \{w_{1,k}, w_{2,k}, w_{3,k}\}, \ldots,  \boldsymbol{\pi}_{k} \, .   
$$    
    
Putting these subsets in $\mathcal{B}_k$, Algorithm \ref{algo:sets}  has to find a new permutation, that is, $(w_{1,k}', \ldots, w_{d_k,k}')$ that avoids some replications described by the first instruction of  Algorithm \ref{algo:sets}. We use the super sets $\mathcal{B}_k$ and $\mathcal{E}_k$ for controlling the choice of new permutations. Formally, consider integers 
\begin{equation} \label{eq:jko}    
j_{0, k} = \left\{  
\begin{array}{cl}
\frac{d_k}{2} & \mbox{if}\; d_k \, \mbox{is even} \\ 
\frac{d_k+1}{2} & \mbox{otherwise} \\  
\end{array}   
\right., \, k=2,\ldots, K \, ,     
\end{equation}     
and the set $\mathcal{A}_{j_{0, k}}$ given by 
$$
\mathcal{A}_{j_{0, k}} =\{ u_k \subseteq \boldsymbol{\pi}_{k}\, : \, |u_k|= j_{0, k} \}, \; k=2,\ldots, K \, .     
$$        
The set $\mathcal{A}_{j_{0, k}}$ is consisted of all the subsets of $\boldsymbol{\pi}_{k}$ that contain exactly $j_{0, k}$ elements, and we can see that the cardinal of $\mathcal{A}_{j_{0, k}}$ is $ \left|\mathcal{A}_{j_{0, k}}\right| = \binom{d_k}{j_{0, k}}$.  The algorithm takes $\mathcal{A}_{j_{0, k}}$ as input and provides the super-sets $\mathcal{B}_k, \mathcal{P}_k, \mathcal{E}_k$. For the first step (i.e., $e_0=1$), Algorithm \ref{algo:sets} focuses on selecting $d_k=\binom{d_k}{e_0}$ permutations with different sets of the form $\{w_{1,k}, \ldots, w_{e_0,k}\}, \ldots, \{w_{1,k}, \ldots, w_{j_{0,k},k}\}$ in one hand and $\{w_{_{d_k-J_{0,k}+1,k}}, \ldots, w_{_{d_k, k}}\}, \ldots, \{w_{_{d_k-e_0+1, k}}, \ldots, w_{_{d_k, k}}\}$ in the other hand. We repeat that process by increasing $e_0$ until we are able to derive all the subsets of $\{1, \ldots, d\}$ using the permutations selected by  Algorithm \ref{algo:sets}, that is, until $\left|\mathcal{A}_{j_{0, k}}\right| =0$. We put the  selected permutations in $\mathcal{P}_k$ (see Lemma \ref{lem:recsets} for more details). The formal algorithm is given as follows:\\    
         
\begin{algorithm}[H]  \label{algo:sets}     
\SetAlgoLined  
 initialization : $\mathcal{B}_k \gets \mathcal{P}_k \gets \mathcal{E}_k \gets \emptyset$; $i \gets  e_0 \gets 1$\;   
 $\mathcal{A}_{j_{0, k}} \gets \left\{u_k \subseteq \boldsymbol{\pi}_{k} : |u_k|=j_{0, k}\right\}$\;
 \While{ $ \left| \mathcal{A}_{j_{0, k}} \right| > 0$}{ 
  Find a permutation $\bo{w}_k=(w_{1,k}, \ldots, w_{d_k,k})$ of $\boldsymbol{\pi}_{k}$ such that : 
	$\left[\{w_{1,k}, \ldots, w_{\jmath, k}  \} \notin \mathcal{B}_k,\, \, \, \jmath =e_0,\ldots, j_{0, k} \right]$ and $\left[\{w_{_{d_k-\jmath+1,k}}, \ldots, w_{_{d_k, k}} \} \notin \mathcal{E},\, \, \, \jmath =e_0,\ldots, j_{0, k} \right]$ and 
	$\left\{w_{1,k}, \ldots, w_{j_{0, k}, k} \right\} \in \mathcal{A}_{j_{0, k}}$\; 
	$\mathcal{A}_{j_{0, k}}  \gets \mathcal{A}_{j_{0, k}}  \setminus \left\{w_{1,k}, \ldots, w_{j_{0, k}, k} \right\}$\; 
		 $\mathcal{B}_k \gets \mathcal{B}_k \bigcup \left\{ \left\{w_{1,k}, \ldots, w_{\jmath, k} \right\},\,  \jmath = e_0,\ldots, d_k -e_0 +1 \right\}$\; 
		$\mathcal{E}_k \gets \mathcal{E}_k \bigcup \left\{ \{w_{\jmath, k}, \ldots, w_{_{d_k,k}}\},\, \jmath =  j_{0, k} +1, \ldots, d_k-e_0 +1 \right\}$\;  
		$\mathcal{P}_k\gets \mathcal{P}_k \bigcup \bo{w}_k$\; 
		$i \gets  i+ 1$\;   
	 \If{$ \binom{d_k}{e_0} <i \leq \binom{d_k}{e_0 + 1}$}{  
			$e_0 \gets e_0 +1 $\;} 
  }   
 \caption{Construction of the sets $\mathcal{B}_k$ and $\mathcal{P}_k$ for all $k \in \{2, \ldots, K\}$.} 
\end{algorithm} 
\vspace{0.5cm}           
 
The set $\mathcal{P}_k$ of Algorithm \ref{algo:sets} contains $\binom{d_k}{j_{0, k}}$ permutations. The set $\mathcal{B}_k$ is built using $\mathcal{P}_k$, and it is consisted of sets containing the first $\imath$ elements of $\bo{w}_k$ with $\imath =1, \ldots, d_k$ and for all $ \bo{w}_k \in \mathcal{P}_k$. Lemma \ref{lem:recsets} gives the properties of the sets $\mathcal{P}_k$ and $\mathcal{B}_k$. 
\begin{lemma} \label{lem:recsets}
Consider $j_{0, k}$ given by (\ref{eq:jko}) and $\mathcal{B}_k,\, \mathcal{P}_k$ given by Algorithm \ref{algo:sets}. Then, 
\begin{eqnarray}  \label{eq:supsetpart}   
  \mathcal{B}_k= \left\{ u_k \subseteq \boldsymbol{\pi}_{k}  : |u_k|>0 \right\} \, ; 
\end{eqnarray}              
\begin{eqnarray}  \label{eq:supsetpart1}    
  \mathcal{B}_k= \left\{ \{w_{1, k},  \ldots, w_{\imath, k} \} ,\, \imath =1,\ldots, d_k : \forall\, \bo{w}_k \in \mathcal{P}_k \right\} \, .    
\end{eqnarray}     
\end{lemma} 
\begin{preuve}         
See Appendix \ref{app:lem:recsets}.  
\begin{flushright}                   
$\Box$  
\end{flushright}
\end{preuve}  

Using Lemma \ref{lem:recsets}, we quantify the necessary and sufficient number of equivalent representations of $\M(\bo{\X})$ given $\bo{\X}_u$ for all $u \subseteq \{1, \ldots, d\}$ in Theorem \ref{theo:tierep}. For $\bo{w}_k \in \mathcal{P}_k$, recall that $\bo{w}_{\sim 1,k} := (w_{2,k}, \ldots, w_{d_k,k})$ and the cardinal $\left|\mathcal{P}_k\right| =\binom{d_k}{j_{0, k}}$ with $k=2, \ldots, K$.      
   
\begin{theorem}  \label{theo:tierep} 
Under (A1), consider integers $p_2 \leq d_2, \ldots, p_K \leq d_K$. Then, \\
       
$\quad$ (i) the minimum number of equivalent representations of $\M(\bo{\X})$ given $\bo{\X}_u$ for all $u \subseteq \{1, \ldots, d\}$ is
\begin{equation} 
R_{\min} := \prod_{ k =2}^K  \binom{d_k}{j_{0, k}} \, .    
\end{equation}     
         
Such representations are given by       
\begin{eqnarray} \label{eq:genrep}  
\M\left(\bo{\X} \right) &\stackrel{d}{=} &  g_{_{\ell}} \left(\bo{\X}_{\Ori}, \X_{w_{1, 2}},  \bo{Z}_{\bo{w}_{\sim 1, 2}}, \ldots, \X_{w_{1, K}},  \bo{Z}_{\bo{w}_{\sim 1, K}}\right) \, ,      
\end{eqnarray}
where $\ell :=\left(\bo{w}_2, \ldots, \bo{w}_K \right)$ for all $\bo{w}_k \in \mathcal{P}_k$ and $k=2,\,\ldots,\, K$.\\  

$\quad$ (ii) The minimum number of representations of $\M(\bo{\X})$  given $\bo{\X}_{u}$  for all \\ 
$u \in \left\{\{\vv_1, \, u_2, \ldots, u_K\} \, : \,\vv_1 \subseteq \Ori; \, \, u_k \subseteq \boldsymbol{\pi}_{k},\, |u_k| \leq p_k \leq j_{0,k},\, k=2,\ldots, K\right\}$ is 
\begin{equation}
R_{p_2\ldots p_K}  := \max_{2 \leq k \leq K} \binom{d_k}{p_k} \, .      
\end{equation}       
\end{theorem}        
\begin{preuve}     
See Appendix \ref{app:theo:tierep}.
\begin{flushright}       
$\Box$
\end{flushright}   
\end{preuve}      
   
Equation (\ref{eq:genrep}) provides generic and equivalent representations of $\M(\bo{\X})$ that share the same distribution (see Lemma \ref{lem:equivdist}). When a function includes only independent variables, we see that $R_{\min} =1$. Reducing $R_{\min}$ will depend on the analysis of interest. For instance,  one representation of $\M$ is sufficient to determine the distribution of $\M(\bo{\X})$ conditional on $\bo{\X}_u$ for all  
$$
u \in \left\{\{\vv_1, \,w_{1, k}, \ldots, w_{\imath, k}\} : \forall\,  \vv_1 \subseteq \Ori; \,  \, \forall\, \imath \in \{0, \ldots, d_k\}, \, k=2,\,\ldots, K \right\} \, .  
$$  
Likewise, $R_1 =\max(d_2, \ldots, d_K)$ equivalent representations are needed for assessing the effects of $\X_j$ for all  $j \in \{1, \ldots, d\}$. It is worth noting that $R_1$ equivalent representations can lead to assess the effects of other groups of inputs. 
										       		 			
\section{Dependent multivariate sensitivity analysis} \label{sec:msa}      
This section extends dGSIs for models with non-independent variables (\cite{lamboni21}) by considering every subset of inputs. 
Since Equation (\ref{eq:genrep}) provides an equivalent representation of a function $\M : \R^d \to \R^\NN$, that is,               
\begin{eqnarray}   
\M\left(\bo{\X}\right) &\stackrel{d}{=} &  g_{_{\ell}} \left(\bo{\X}_{\Ori}, \X_{w_{1, 2}},  \bo{Z}_{\bo{w}_{\sim 1, 2}}, \ldots, \X_{w_{1, K}},  \bo{Z}_{\bo{w}_{\sim 1, K}}\right) \, ,   \nonumber   
\end{eqnarray} 
and  $g_{_{\ell}}$ includes only independent variables, multivariate sensitivity anlysis (\cite{lamboni11,gamboa14,lamboni18a,lamboni18,perrin21,lamboni22}) can be used for defining GSIs of these variables. Thanks to Lemmas \ref{lem:equivdist}-\ref{lem:recsets}, we are able to define dGSIs of each subset of inputs using multivariate sensitivity anlysis. To ensure the existence of dGSIs, we assume that\\    
       
 (A2): $\displaystyle 0 < \esp\left[\norme{\M\left(\bo{\X}\right)}^2 \right] < + \infty$.\\

Recall that for all $0\leq  p_k \leq d_k$ (see Equation (\ref{eq:setvuk})), 
$$
\rr_k =(w_{2, k}, \ldots, w_{p_k, k}); \quad \vv_k =(w_{1, k}, \rr_k); \quad \boldsymbol{\varpi}_k :=(w_{p_k+1,k}, \ldots, w_{d_k,k}) \, ,  
$$  
 with $k=2, \ldots, K$. According to Lemmas \ref{lem:equivdist}-\ref{lem:recsets}, for all $u \subseteq \{1, \ldots, d\}$, there exists $\vv_1 \subseteq \Ori$, a permutation $\bo{w}_k \in \mathcal{P}_k$ and $p_k$ with $k=2, \ldots, K$ such that 
$$
u = \left\{\vv_1, \vv_2, \ldots, \vv_K  \right\};
\qquad     
\bo{\X}_u = \left(\bo{\X}_{\vv_1}, \, \bo{\X}_{\vv_2},\ldots , \bo{\X}_{\vv_K} \right) \, . 
$$ 
It is worth noting that dGSIs defined and studied in \cite{lamboni21} concern the inputs $\bo{\X}_u$ with $0\leq p_k \leq 1$ for all $k \in\{2,\ldots, K\}$. In this paper, we consider $0\leq p_k \leq d_k$ to cope with every subset of inputs such as two inputs from one block of dependent variables. We then study the properties of such indices and construct their estimations. 
          
\subsection{Extended definition of dependent generalized sensitivity indices} \label{sec:extdef} 
The definitions of GSIs and dGSIs are based on sensitivity functionals (SFs), which contain the primary information about the single and overall contributions of inputs over the whole model outputs (\cite{lamboni16,lamboni16b,lamboni18,lamboni18a,lamboni19,lamboni20,lamboni21}).  
To define SFs of $\bo{\X}_u$, we use $\oo := \{\Ori,\Select \}$ and  $\bo{\X}_{\oo} := (\X_i, \, \forall\, i \in \oo)$; $   
\bo{Z}_{\rr} :=\left(\bo{Z}_{\rr_2}, \ldots, \bo{Z}_{\rr_K} \right)
$; 
 and 
$
\bo{Z}_{\boldsymbol{\varpi}} :=\left(\bo{Z}_{\boldsymbol{\varpi}_2}, \ldots, \bo{Z}_{\boldsymbol{\varpi}_K} \right)
$    
for concise notations. Thus, we can see that $\bo{Z}_{\bo{w}_{\sim 1}} =\left(\bo{Z}_{\bo{w}_{\sim 1, 2}}, \ldots, \bo{Z}_{\bo{w}_{\sim 1, K}} \right) = \left(\bo{Z}_{\rr},  \bo{Z}_{\boldsymbol{\varpi}} \right)$, and Lemma \ref{lem:equivdist} becomes 
     
\begin{equation}   \label{eq:equirep021} 
\M(\bo{\X}) |\, \bo{\X}_u \stackrel{d}{=} g_{_{\ell}} \left(\bo{\X}_{\oo}, \bo{Z}_{\rr},  \bo{Z}_{\boldsymbol{\varpi}} \right) |\, \bo{\X}_{\oo \cap u}, \bo{Z}_{\rr} \, ,   
\end{equation}
because $\oo \cap u =\{\vv_1, \Select \cap u\}$. \\   
	
The first-order SF of $\bo{\X}_u$ with $u=\left\{\vv_1,\, \vv_{2}, \ldots, \vv_{K} \right\}$ is given by 
\begin{equation}     \label{eq:fosf}  
\M^{fo}_u \left(\bo{\X}_{\oo \cap u}, \bo{Z}_{\rr} \right)
 := \esp \left[ g_{_{\ell}} \left(\bo{\X}_{\oo}, \bo{Z}_{\bo{w}_{\sim 1}} \right) |\, \bo{\X}_{\oo \cap u}, \bo{Z}_{\rr} \right] 
- \esp\left[g_{_{\ell}} \left(\bo{\X}_{\oo}, \bo{Z}_{\bo{w}_{\sim 1}}\right)\right] \, ,   
\end{equation}           
and the total SF, which contains the overall information about $\bo{\X}_u$, is given by 
\begin{equation}       \label{eq:totsf}
 \M^{tot}_u\left(\bo{\X}_{\oo}, \bo{Z}_{\bo{w}_{\sim 1}}\right) := g_{_{\ell}} \left(\bo{\X}_{\oo}, \bo{Z}_{\bo{w}_{\sim 1}}\right) 
- \esp_{\bo{\X}_{\oo \cap u}, \bo{Z}_{\rr}} \left[ g_{_{\ell}} \left(\bo{\X}_{\oo}, \bo{Z}_{\bo{w}_{\sim 1}} \right)\right] \, ,  
\end{equation}            
where $\esp_{\bo{\X}_{\oo \cap u}, \bo{Z}_{\rr}}$ means that the expectation is taken w.r.t. $\bo{\X}_{\oo \cap u}, \bo{Z}_{\rr}$. \\ 
                          
The SFs given by (\ref{eq:fosf})-(\ref{eq:totsf}) are random vectors, and their components may be correlated and/or dependent. Using the variance-covariance as importance measure, the covariance  of $\M^{fo}_u \left(\bo{\X}_{\oo \cap u}, \bo{Z}_{\rr}\right)$  also known as the first-order covariance is given by 
\begin{equation}    \label{eq:siguf1}    
\D_u := \esp\left[\M^{fo}_u\left(\bo{\X}_{\oo \cap u}, \bo{Z}_{\rr} \right) \M^{fo}_u\left(\bo{\X}_{\oo \cap u}, \bo{Z}_{\rr}\right)^\T \right] \, . \nonumber  
\end{equation}                                    
Likewise, the covariances of $\M^{tot}_u\left(\bo{\X}_{\oo}, \bo{Z}_{\bo{w}_{\sim 1}}\right)$ and the model outputs are given as follows:   
\begin{equation}    \label{eq:sigut1}  
 \D_u^{tot}  :=  \esp\left[\M^{tot}_u\left(\bo{\X}_{\oo}, \bo{Z}_{\bo{w}_{\sim 1}}\right) 
 \M^{tot}_u\left(\bo{\X}_{\oo}, \bo{Z}_{\bo{w}_{\sim 1}}\right)^\T \right] \, ,  \nonumber 
\end{equation}  
\begin{eqnarray}  \label{eq:siguti1}     
 \Sigma &:=&  \esp\left[g_{_{\ell}} \left(\bo{\X}_{\oo}, \bo{Z}_{\bo{w}_{\sim 1}} \right)
 g_{_{\ell}} \left(\bo{\X}_{\oo}, \bo{Z}_{\bo{w}_{\sim 1}}\right)^\T \right]
 - \esp\left[g_{_{\ell}} \left(\bo{\X}_{\oo}, \bo{Z}_{\bo{w}_{\sim 1}} \right) \right]
	\esp\left[g_{_{\ell}} \left(\bo{\X}_{\oo}, \bo{Z}_{\bo{w}_{\sim 1}}\right)^\T \right]  \, .  \nonumber     
\end{eqnarray}        
   
Based on the above covariances, Definition \ref{def:gsi} extends dGSIs introduced in \cite{lamboni21}.

\begin{defi}    \label{def:gsi} 
Consider $\bo{\X}_u$ with $u=\left\{\vv_1,\, \vv_{2}, \ldots, \vv_{K} \right\}$, and assume (A2) holds.\\ 
   
$\quad$ (i) The first-type dGSIs of $\bo{\X}_u$ are given as follows:

\begin{equation}        \label{eq:gsif10}   
dGSI_u^{1, M} := \frac{\trace\left(\D_u \right)}{\trace\left(\Sigma\right)}  
\, , 
\qquad \quad  
dGSI_{T_u}^{1, M} := \frac{\trace\left(\D_u^{tot} \right)}{\trace\left(\Sigma\right)} \, .   
\end{equation}  
          
$\quad$ (ii) The second-type dGSIs  of $\bo{\X}_u$ are defined as follows: 
\begin{equation}        \label{eq:gsif11}      
dGSI_u^{2, M} := \frac{\normf{\D_u}}{\normf{\Sigma}} \, ,
\qquad \quad 
dGSI_{T_u}^{2, M} := \frac{\normf{\D_u^{tot}}}{\normf{\Sigma}} \, .
\end{equation}          
\end{defi}  
   
Note that $dGSI_u^{1, M}$ (resp. $dGSI_u^{2, M}$) is the first-order index of $\bo{\X}_u$;  $dGSI_{T_u}^{1, M}$ (resp. $dGSI_{T_u}^{2, M}$) is the total index of $\bo{\X}_u$; and the second-type dGSIs account for the correlations among the components of SFs. When $\NN =1$, both types of dGSIs  are equal and boil down to  dependent sensitivity indices (dSIs) (see Defintion \ref{def:dsi}).
       
\begin{defi}  \label{def:dsi} 
For a function $\M :\R^d \to \R$ ($\NN =1$), the dSIs of $\bo{X}_u$ are given by   
\begin{equation}        \label{eq:gsif10}       
dS_u := \frac{\D_u}{\Sigma} \, ,   
\qquad \quad 
dS_{T_u} := \frac{\D_u^{tot}}{\Sigma} \, .   
\end{equation}     
\end{defi}    
          
\subsection{Properties of dependent generalized sensitivity indices} \label{sec:prop} 
Both types of dGSIs share the same properties as those proposed in \cite{lamboni21} for all $0 \leq p_k \leq 1$. Proposition \ref{prop:gsipro} extends such properties for all $p_k \in \{0, \ldots, d_k\}$.  
\begin{prop} \label{prop:gsipro}
Consider an orthogonal matrix $\mathcal{V} \in \R^{\NN \times \NN}$ (i.e., $\mathcal{V}\mathcal{V}^\T=\mathcal{V}^\T\mathcal{V}=\mathcal{I}$), and assume (A2) holds. Then, we have
   \begin{eqnarray}
	     0  \leq  dGSI_{u}^{1, M} \leq  dGSI_{T_u}^{1, M}  \leq  1 \, ,  \\
		   0  \leq  dGSI_{u}^{2, M}  \leq  dGSI_{T_u}^{2, M} \leq  1  \, ; 
	 \end{eqnarray}
  \begin{eqnarray}    
	    dGSI_{u}^{i, M} (\mathcal{V}\M) &=& dGSI_{u}^{i, M} \, , \\
		  dGSI_{T_u}^{i, M} (\mathcal{V}\M) &=& dGSI_{T_u}^{i, M}  \, ,
	\end{eqnarray}   
	where $i=1, \, 2$ denotes the first-type and second-type dGSIs, respectively. 
\end{prop} 
\begin{preuve}
See Appendix  \ref{app:prop:gsipro}.  
\begin{flushright}         
$\Box$  
\end{flushright}                  
\end{preuve} 
   
When the total dGSI of $\bo{\X}_{u}$ is zero or almost zero, we must fix $\bo{\X}_{u}$ using the dependency functions of each $\bo{\X}_{\boldsymbol{\pi}_k}$ with  $k=2, \ldots, K$. Indeed, since we can write $\bo{\X}_{u_k} = r_{\sim u_k}(\bo{\X}_{\sim u_k},\, \bo{Z}_{u_k})$, fixing $\bo{\X}_{u_K}$ comes down to fix $\bo{Z}_{u_k}$ to its nominal values. Instead of fixing $\bo{Z}_{u_k}$, one may also take the expectation over $\bo{Z}_{u_k}$. Moreover, as one can compute the total $dGSI$ of $\bo{\X}_{\boldsymbol{\pi}_k}$ for all $k\geq 2$ using any equivalent representations of $\M$, it becomes possible to quickly identify the non-influential blocks of dependent variables, and then put our computational efforts on the most important groups of dependent variables. \\
     
Regarding the ranking of input variables, Proposition \ref{prop:equdgsi} provides conditions that guarantee equivalent ranking of inputs using either $dGSI_{T_u}^{1, M}$ or $dGSI_{T_u}^{2, M}$. To that end, we use $\mathcal{A}_1 \matleq \mathcal{A}_2$ to say that $\mathcal{A}_2 -\mathcal{A}_1$ is positive semi-definite, also known as the Loewner partial ordering between matrices. 
  
\begin{prop}  \label{prop:equdgsi}  Let $\bo{\X}_u$ and $\bo{\X}_{\omega}$ be two subsets of $\bo{\X}$ having $\D_u^{tot}$, $\D_\omega^{tot}$ as the total-effect covariances.\\ 
  
If $\D_u^{tot} \matleq \D_\omega^{tot}$, then we have
\begin{equation}
dGSI_{T_u}^{1, M} \leq  dGSI_{T_\omega}^{1, M}, \qquad\,  dGSI_{T_u}^{2, M} \leq dGSI_{T_\omega}^{2, M} \, . \nonumber 
\end{equation}  
\end{prop} 
   
When assumption $\D_u^{tot} \matleq \D_\omega^{tot}$ is not satisfied, we may have different ranking of inputs using both types of dGSIs (see Section 6.3 in \cite{lamboni21}). 

\begin{rem}
The definition of dGSIs for he multivariate and functional outputs provided in \cite{lamboni21} (i.e., $\M(\bo{\X}, \theta) \in \R^\NN$ with $\theta \in \Theta \subseteq \R$) can also be extended to cope with every subsets of inputs.  
\end{rem}     
 
\subsection{Estimators of dependent generalized sensitivity indices} \label{sec:est}
In this section, we provide minimum variance and unbiased estimators of the first-order, total-effect covariances and the covariance of the model outputs; consistent estimators of dGSIs of each subset of inputs and their asymptotic distributions. \\ 
 
To derive new expressions of the first-order and total-effect covariances associated with 
$\bo{\X}_u$ that are based on $\M$ (see Proposition \ref{prop:expcovkn}), let us recall that $\bo{Z}_{\bo{w}_{\sim 1}}=\left(\bo{Z}_{\rr},  \bo{Z}_{\boldsymbol{\varpi}} \right)$,
$$
 g_{_{\ell}} \left(\bo{\X}_{\oo},  \bo{Z}_{\bo{w}_{\sim 1}} \right) =   \M\left(\bo{\X}_{\oo},  r_\Select \left(\bo{\X}_\Select,  \bo{Z}_{\bo{w}_{\sim 1}} \right) \right) =   \M\left(\bo{\X}_{\oo},  r_\Select \left(\bo{\X}_\Select, \bo{Z}_{\rr},  \bo{Z}_{\boldsymbol{\varpi}} \right) \right)  \, . 
$$ 
 
\begin{prop} \label{prop:expcovkn}    
Let $\left(\bo{\X}_{\oo}^{(1)},\,  \bo{Z}_{\bo{w}_{\sim 1}}^{(1)} \right)$ and $\left(\bo{\X}_{\oo}^{(2)},\,  \bo{Z}_{\bo{w}_{\sim 1}}^{(2)} \right)$ be two i.i.d. copies of $\left(\bo{\X}_{\oo},\,  \bo{Z}_{\bo{w}_{\sim 1}} \right)$ and $\bo{\X}^{(\imath)} :=\left(\bo{\X}_{\oo}^{(\imath)},\, r_{\Select}\left(\bo{\X}_{\Select}^{(\imath)},\, \bo{Z}_{\bo{w}_{\sim 1}}^{(\imath)} \right) \right)$ with $\imath=1, 2$. Assume  (A2) holds.\\ 
   
$\quad$ (i)  The first-order covariance is given by 
\begin{equation}    \label{eq:siguf1kn}      
\displaystyle  
\D_u =
 \esp\left[\M\left(\bo{\X}^{(1)}\right) \M\left(\bo{\X}_{\oo\cap u}^{(1)},\, \bo{\X}_{\oo\setminus u}^{(2)},\, r_{\Select}\left(\bo{\X}_{\Select \cap u}^{(1)}, \bo{\X}_{\Select \setminus u}^{(2)}, \bo{Z}_{\rr}^{(1)}, \bo{Z}_{\boldsymbol{\varpi}}^{(2)} \right) \right)^\T \right] - \esp\left[\M(\bo{\X}^{(1)}) \right] \esp\left[\M(\bo{\X}^{(1)}) \right]^\T  \, .      \nonumber           
\end{equation}                     
            
$\quad$ (ii) The total-effect covariance is given by     
\begin{equation}    \label{eq:sigut1kn} 
\displaystyle         
\D_u^{tot} = \esp\left[\M(\bo{\X}^{(1)}) \M(\bo{\X}^{(1)})^\T\right] -
\esp\left[\M\left(\bo{\X}^{(1)}\right) \M\left(\bo{\X}_{\oo\cap u}^{(2)},\, \bo{\X}_{\oo\setminus u}^{(1)},\, r_{\Select}\left(\bo{\X}_{\Select \cap u}^{(2)}, \bo{\X}_{\Select \setminus u}^{(1)}, \bo{Z}_{\rr}^{(2)}, \bo{Z}_{\boldsymbol{\varpi}}^{(1)} \right) \right)^\T \right]  \, .      \nonumber             
\end{equation}            
\end{prop}      
\begin{preuve} 
Since the components of $\bo{\X}_{\oo}, \bo{Z}_{\rr},  \bo{Z}_{\boldsymbol{\varpi}}$ are independent and bearing in mind Equation (\ref{eq:equirep021}), the proof is similar to the proof of Proposition 3 provided in \cite{lamboni21}.    
\begin{flushright}
$\Box$              
\end{flushright}  
\end{preuve}
  
Based on Proposition \ref{prop:expcovkn}, we derive the  minimum variance and unbiased estimators (MVUEs) of the variance-covariances of SFs in Theorem \ref{theo:sigestM}. To provide such results, we are given two i.i.d. samples, that is,  $\left\{\left(\bo{\X}_{i, \oo}^{(1)},\,  \bo{Z}_{i, \bo{w}_{\sim 1}}^{(1)} \right) \right\}_{i=1}^m$ and $\left\{\left(\bo{\X}_{i, \oo}^{(2)},  \bo{Z}_{i, \bo{w}_{\sim 1}}^{(2)} \right) \right\}_{i=1}^m$  from $\left(\bo{\X}_{\oo},\,  \bo{Z}_{\bo{w}_{\sim 1}} \right)$, and we use $\bo{\X}^{(\jmath)}_i :=\left(\bo{\X}_{i, \oo}^{(\jmath)},\, r_{\Select}\left(\bo{\X}_{i,\Select}^{(\jmath)},\, \bo{Z}_{i, \bo{w}_{\sim 1}}^{(\jmath)} \right) \right)$ with $\jmath=1, 2$, $i=1, \ldots, m$. Moreover, we use $\mathcal{X}_i :=\left(\bo{\X}_{i, \oo}^{(1)},\, \bo{\X}_{i, \oo}^{(2)},\, \bo{Z}_{i, \bo{w}_{\sim 1}}^{(1)}, \, \bo{Z}_{i, \bo{w}_{\sim 1}}^{(2)} \right)$ with $i=1,\ldots, m$, and we consider the following symmetric kernels:
\begin{eqnarray} 
&& K(\mathcal{X}_i) := \left(
 \left[\M\left(\bo{\X}_{i}^{(1)}\right) - \M\left(\bo{\X}_{i, \oo\cap u}^{(2)},\, \bo{\X}_{i, \oo\setminus u}^{(1)},\, r_{\Select}\left(\bo{\X}_{i, \Select \cap u}^{(2)}, \bo{\X}_{i, \Select \setminus u}^{(1)}, \,
\bo{Z}_{i, \rr}^{(2)},\, \bo{Z}_{i, \boldsymbol{\varpi}}^{(1)} \right) \right) \right] \right.  \nonumber \\ 
& &   \times \left[\M\left(\bo{\X}_{i, \oo\cap u}^{(1)},\, \bo{\X}_{i, \oo\setminus u}^{(2)},\, r_{\Select}\left(\bo{\X}_{i, \Select \cap u}^{(1)}, \bo{\X}_{i, \Select \setminus u}^{(2)}, \,   
\bo{Z}_{i, \rr}^{(1)},\, \bo{Z}_{i, \boldsymbol{\varpi}}^{(2)} \right) \right) - \M\left(\bo{\X}_{i}^{(2)} \right)\right]^\T \nonumber \\    
& & +  \left[\M\left(\bo{\X}_{i, \oo\cap u}^{(1)},\, \bo{\X}_{i, \oo\setminus u}^{(2)},\, r_{\Select}\left(\bo{\X}_{i, \Select \cap u}^{(1)}, \bo{\X}_{i, \Select \setminus u}^{(2)}, \,   
\bo{Z}_{i, \rr}^{(1)},\, \bo{Z}_{i, \boldsymbol{\varpi}}^{(2)} \right) \right) - \M\left(\bo{\X}_{i}^{(2)} \right)\right]    \nonumber \\
& & \left. \times   \left[\M\left(\bo{\X}_{i}^{(1)}\right) - \M\left(\bo{\X}_{i, \oo\cap u}^{(2)},\, \bo{\X}_{i, \oo\setminus u}^{(1)},\, r_{\Select}\left(\bo{\X}_{i, \Select \cap u}^{(2)}, \bo{\X}_{i, \Select \setminus u}^{(1)}, \,
\bo{Z}_{i, \rr}^{(2)},\, \bo{Z}_{i, \boldsymbol{\varpi}}^{(1)} \right) \right) \right]^\T  
 \right) \, , \nonumber         
\end{eqnarray}   
\begin{eqnarray} 
& & K^{tot}(\mathcal{X}_i) := \left( \left[\M\left(\bo{\X}_{i}^{(1)}\right) - \M\left(\bo{\X}_{i, \oo\cap u}^{(2)},\, \bo{\X}_{i, \oo\setminus u}^{(1)},\, r_{\Select}\left(\bo{\X}_{i, \Select \cap u}^{(2)}, \bo{\X}_{i, \Select \setminus u}^{(1)}, \,
\bo{Z}_{i, \rr}^{(2)},\, \bo{Z}_{i, \boldsymbol{\varpi}}^{(1)} \right) \right) \right]  \right. \nonumber \\  
& &  \times \left[\M\left(\bo{\X}_{i}^{(1)}\right) - \M\left(\bo{\X}_{i, \oo\cap u}^{(2)},\, \bo{\X}_{i, \oo\setminus u}^{(1)},\, r_{\Select}\left(\bo{\X}_{i, \Select \cap u}^{(2)}, \bo{\X}_{i, \Select \setminus u}^{(1)}, \,
\bo{Z}_{i, \rr}^{(2)},\, \bo{Z}_{i, \boldsymbol{\varpi}}^{(1)} \right) \right) \right]^\T  \nonumber \\     
& & +  \left[\M\left(\bo{\X}_{i, \oo\cap u}^{(1)},\, \bo{\X}_{i, \oo\setminus u}^{(2)},\, r_{\Select}\left(\bo{\X}_{i, \Select \cap u}^{(1)}, \bo{\X}_{i, \Select \setminus u}^{(2)}, \,    
\bo{Z}_{i, \rr}^{(1)},\, \bo{Z}_{i, \boldsymbol{\varpi}}^{(2)} \right) \right) - \M\left(\bo{\X}_{i}^{(2)} \right)\right]   \nonumber \\
& &   \times \left. \left[\M\left(\bo{\X}_{i, \oo\cap u}^{(1)},\, \bo{\X}_{i, \oo\setminus u}^{(2)},\, r_{\Select}\left(\bo{\X}_{i, \Select \cap u}^{(1)}, \bo{\X}_{i, \Select \setminus u}^{(2)}, \,   
\bo{Z}_{i, \rr}^{(1)},\, \bo{Z}_{i, \boldsymbol{\varpi}}^{(2)} \right) \right) - \M\left(\bo{\X}_{i}^{(2)} \right)\right]^\T \right)   \, .  \nonumber          
\end{eqnarray}           
                      
\begin{theorem}        \label{theo:sigestM}  
Assume (A2) holds and $\M(\bo{\X})$ has finite fourth moments (A3). \\
 
 $\quad (i) $ The consistent and MVUE of $\D_u$ is given by
\begin{equation} \label{eq:sigfoestM} 
  \widehat{\D_u}  :=  \frac{1}{4m}\sum_{i=1}^m K(\mathcal{X}_i) \xrightarrow{\mathcal{P}} \D_u \, ,
\end{equation}   
where $\xrightarrow{\mathcal{P}}$ stands for the convergence in probability.\\
                 
$\quad (ii) $ The consistent and MVUE of $\D_u^{tot}$ is given by               
\begin{equation} \label{eq:sigtotestM}
  \widehat{\D_u^{tot}}  := \frac{1}{4m}\sum_{i=1}^m  K^{tot}(\mathcal{X}_i) \xrightarrow{\mathcal{P}}  \D_u^{tot} \, .            
\end{equation}                                      

$\quad (iii) $ The consistent and MVUE of $\Sigma$ is given by  
\begin{eqnarray} \label{eq:sigvestM}
 \widehat{\Sigma}  &:=& \frac{1}{2m}\sum_{i=1}^m  
 \left[\M\left(\bo{\X}_{i}^{(1)}\right) - \M\left(\bo{\X}_{i}^{(2)} \right) \right]  \left[\M\left(\bo{\X}_{i}^{(1)} \right) - \M\left(\bo{\X}_{i}^{(2)}\right) \right]^\T  \xrightarrow{\mathcal{P}}  \Sigma \, .   
\end{eqnarray}                  
\end{theorem}  
\begin{preuve} 
See Appendix \ref{app:theo:sigestM}.  
\begin{flushright}  
$\Box$       
\end{flushright}     
\end{preuve}    

For single response models (i.e., $\NN=1$), the consistent and MVUEs of the covariances of SFs provided in Theorem \ref{theo:sigestM} have simple expressions (see Corollary \ref{coro:rem:sigsis}).  
\begin{corollary} \label{coro:rem:sigsis}  
Assume (A2)-(A3) hold and  $\NN=1$. Then, the MVUEs $\widehat{\D_u}$, $\widehat{\D_u^{tot}}$ and   $\widehat{\Sigma}$ boil down to 

\begin{eqnarray} \label{eq:sigfoests}
 \widehat{\sigma_u}  &:=&  \frac{1}{2m}\sum_{i=1}^m  \left[\M\left(\bo{\X}_{i}^{(1)}\right) - \M\left(\bo{\X}_{i, \oo\cap u}^{(2)},\, \bo{\X}_{i, \oo\setminus u}^{(1)},\, r_{\Select}\left(\bo{\X}_{i, \Select \cap u}^{(2)}, \bo{\X}_{i, \Select \setminus u}^{(1)}, \,
\bo{Z}_{i, \rr}^{(2)},\, \bo{Z}_{i, \boldsymbol{\varpi}}^{(1)} \right) \right) \right]  \nonumber \\  
& & \times \left[\M\left(\bo{\X}_{i, \oo\cap u}^{(1)},\, \bo{\X}_{i, \oo\setminus u}^{(2)},\, r_{\Select}\left(\bo{\X}_{i, \Select \cap u}^{(1)}, \bo{\X}_{i, \Select \setminus u}^{(2)}, \,   
\bo{Z}_{i, \rr}^{(1)},\, \bo{Z}_{i, \boldsymbol{\varpi}}^{(2)} \right) \right) - \M\left(\bo{\X}_{i}^{(2)} \right)\right]  \, ;   \nonumber  
\end{eqnarray}              
\begin{eqnarray} \label{eq:sigtotests}
 \widehat{\sigma_u^{tot}}  &:=& \frac{1}{4m}\sum_{i=1}^m  \left( \left[\M\left(\bo{\X}_{i}^{(1)}\right) - \M\left(\bo{\X}_{i, \oo\cap u}^{(2)},\, \bo{\X}_{i, \oo\setminus u}^{(1)},\, r_{\Select}\left(\bo{\X}_{i, \Select \cap u}^{(2)}, \bo{\X}_{i, \Select \setminus u}^{(1)}, \,
\bo{Z}_{i, \rr}^{(2)},\, \bo{Z}_{i, \boldsymbol{\varpi}}^{(1)} \right) \right) \right]^2 \right. \nonumber \\  
& &   +  \left. 
\left[\M\left(\bo{\X}_{i, \oo\cap u}^{(1)},\, \bo{\X}_{i, \oo\setminus u}^{(2)},\, r_{\Select}\left(\bo{\X}_{i, \Select \cap u}^{(1)}, \bo{\X}_{i, \Select \setminus u}^{(2)}, \,   
\bo{Z}_{i, \rr}^{(1)},\, \bo{Z}_{i, \boldsymbol{\varpi}}^{(2)} \right) \right) - \M\left(\bo{\X}_{i}^{(2)} \right)\right]^2 \right)  \, ;  \nonumber
\end{eqnarray}                                 
\begin{equation} \label{eq:sigvests}
 \widehat{\sigma}  := \frac{1}{2m}\sum_{i=1}^m  
 \left[\M\left(\bo{\X}_{i}^{(1)}\right) - \M\left(\bo{\X}_{i}^{(2)} \right) \right]^2 \, .  \nonumber 
\end{equation}                 
\end{corollary} 
 \begin{preuve} 
The proof is obvious using Theorem \ref{theo:sigestM}. 
 \begin{flushright} 
 $\Box$    
 \end{flushright}         
 \end{preuve} 

\begin{rem}
When one is only interested in the total-effect covariances, the following unbiased and consistent estimators of $\Sigma_u^{tot}$ are sufficient:
\begin{eqnarray} 
\widehat{\Sigma_u^{tot '}} &:= & \frac{1}{2m}\sum_{i=1}^m 
 \left[\M\left(\bo{\X}_{i}^{(1)}\right) - \M\left(\bo{\X}_{i, \oo\cap u}^{(2)},\, \bo{\X}_{i, \oo\setminus u}^{(1)},\, r_{\Select}\left(\bo{\X}_{i, \Select \cap u}^{(2)}, \bo{\X}_{i, \Select \setminus u}^{(1)}, \, 
\bo{Z}_{i, \rr}^{(2)},\, \bo{Z}_{i, \boldsymbol{\varpi}}^{(1)} \right) \right) \right]  \nonumber \\  
& &  \times \left[\M\left(\bo{\X}_{i}^{(1)}\right) - \M\left(\bo{\X}_{i, \oo\cap u}^{(2)},\, \bo{\X}_{i, \oo\setminus u}^{(1)},\, r_{\Select}\left(\bo{\X}_{i, \Select \cap u}^{(2)}, \bo{\X}_{i, \Select \setminus u}^{(1)}, \,
\bo{Z}_{i, \rr}^{(2)},\, \bo{Z}_{i, \boldsymbol{\varpi}}^{(1)} \right) \right) \right]^\T  \, ;  \nonumber   
\end{eqnarray} 
\begin{equation}   
\widehat{\sigma_u^{tot '}} := \frac{1}{2m}\sum_{i=1}^m 
\left[\M\left(\bo{\X}_{i}^{(1)}\right) - \M\left(\bo{\X}_{i, \oo\cap u}^{(2)},\, \bo{\X}_{i, \oo\setminus u}^{(1)},\, r_{\Select}\left(\bo{\X}_{i, \Select \cap u}^{(2)}, \bo{\X}_{i, \Select \setminus u}^{(1)}, \,
\bo{Z}_{i, \rr}^{(2)},\, \bo{Z}_{i, \boldsymbol{\varpi}}^{(1)} \right) \right) \right]^2  \, .  \nonumber  
\end{equation}
\end{rem}
        
For computing the first-order and total-effect covariances related to $\bo{\X}_u$ for all $u \subseteq \{1, \ldots, d\}$, different model runs are needed. Some of such model runs can be combined and then used for computing the model outputs covariance. In what follows, we use a sample of size $M \geq m$  for computing the outputs covariance given by (\ref{eq:sigvestM}). The operator $\vect(\cdot)$ transforms a matrix $\Sigma \in \R^{\NN \times \NN}$ into a vector, that is, $\vect(\Sigma) \in \R^{\NN^2}$ and  $\mathsf{O} \in \R^{\NN \times \NN}$ denotes the null matrix. Theorem \ref{theo:estgsiM} and Corollary \ref{coro:siests} provide the estimators of dGSIs and dSIs, respectively. 
  
\begin{theorem} \label{theo:estgsiM} 
Assume that (A2)-(A3) hold, $m \to +\infty$,  $M \to +\infty$.\\

$\quad$ (i) The estimators of the first-type dGSIs are given  as follows:
  
\begin{equation}        
\widehat{dGSI_u^{1, M}} := \frac{\trace\left(\widehat{\D_u}\right)}{\trace\left(\widehat{\Sigma}\right)} \xrightarrow{\mathcal{P}}  dGSI_u^{1, M}  \, ;
\end{equation}      
\begin{equation}
\widehat{dGSI_{T_u}^{1, M}} := \frac{\trace\left(\widehat{\D_u^{tot}}\right)}{\trace\left(\widehat{\Sigma}\right)} \xrightarrow{\mathcal{P}}   dGSI_{T_u}^{1, M} \, .
\end{equation}   
If $m/M  \to 0$, then we have the following asymptotic distributions 
\begin{equation}
\sqrt{m}\left(\widehat{dGSI_u^{1, M}} - dGSI_u^{1, M} \right) \xrightarrow{\mathcal{D}}  \mathcal{N} \left(0, \frac{\var\left[\trace(K(\mathcal{X}_1))\right]}{(\trace(\Sigma))^2} \right) \, ;   \nonumber   
\end{equation}  
\begin{equation}
\sqrt{m}\left(\widehat{dGSI_{T_u}^{1, M}} - dGSI_{T_u}^{1, M} \right) \xrightarrow{\mathcal{D}}  \mathcal{N} \left(0, \frac{\var\left[\trace(K^{tot}(\mathcal{X}_1))\right]}{(\trace(\Sigma))^2} \right) \, .  \nonumber 
\end{equation}    
$\quad$ (ii)  For the second-type dGSIs, we have 
\begin{equation}         
\widehat{dGSI_u^{2, M}} := \frac{\normf{\widehat{\D_u}}}{\normf{\widehat{\Sigma}}}
\xrightarrow{\mathcal{P}} dGSI_u^{2, M} \, ,  
\end{equation}
\begin{equation}
\widehat{dGSI_{T_u}^{2, M}} := \frac{\normf{\widehat{\D_u^{tot}}}}{\normf{\widehat{\Sigma}}}
\xrightarrow{\mathcal{P}} dGSI_{T_u}^{2, M} \, .    
\end{equation}  
If $m/M  \to 0$, $\D_u \neq \mathsf{O}$ and $\D_u^{tot} \neq \mathsf{O}$, then we have
\begin{equation}
\sqrt{m}\left(\widehat{dGSI_{u}^{2, M}} - dGSI_{u}^{2, M} \right) \xrightarrow{\mathcal{D}}  \mathcal{N} \left(0, \frac{\vect(\D_u)^\T \var\left[\vect(K(\mathcal{X}_1))\right] \vect(\D_u)}{\normf{\D_u}^2 \normf{\Sigma}^2} \right) \, .   \nonumber
\end{equation}          
\begin{equation}
\sqrt{m}\left(\widehat{dGSI_{T_u}^{2, M}} - dGSI_{T_u}^{2, M} \right) \xrightarrow{\mathcal{D}}  \mathcal{N} \left(0, \frac{\vect(\D_u^{tot})^\T \var\left[\vect(K^{tot}(\mathcal{X}_1))\right] \vect(\D_u^{tot})}{\normf{\D_u^{tot}}^2 \normf{\Sigma}^2} \right) \, .  \nonumber   
\end{equation}        
\end{theorem}  
\begin{preuve} 
See Appendix \ref{app:theo:estgsiM}.   
\begin{flushright}
$\Box$
\end{flushright} 
\end{preuve}   
   
For real-valued functions, the estimators of dSIs (see Definition \ref{def:dsi}) are derived in Corollary~\ref{coro:siests} using Theorem \ref{theo:estgsiM} and Corollary \ref{coro:rem:sigsis}. 
\begin{corollary} \label{coro:siests}   
Assume that (A2)-(A3) hold,  $\NN=1$, $m, M \to +\infty$. The estimators of dSIs are given as follows:   
\begin{equation}   
\widehat{dS_u} := \frac{\widehat{\sigma_u}}{\widehat{\sigma}} \xrightarrow{\mathcal{P}} dS_u; \qquad  \qquad   
\widehat{dS_{T_u}} := \frac{\widehat{\sigma_u^{tot}}}{\widehat{\sigma}} \xrightarrow{\mathcal{P}} dS_{T_u}
\, . 
\end{equation} 
If $m/M  \to 0$, then we have 
\begin{equation}
\sqrt{m}\left(\widehat{dS_u} - dS_u \right) \xrightarrow{\mathcal{D}}  \mathcal{N} \left(0, \frac{\var\left[K(\mathcal{X}_1)\right]}{\sigma^4} \right); 
\quad 
\sqrt{m}\left(\widehat{dS_{T_u}} - dS_{T_u} \right) \xrightarrow{\mathcal{D}}  \mathcal{N} \left(0, \frac{\var\left[K^{tot}(\mathcal{X}_1)\right]}{\sigma^4} \right)
\, .   \nonumber 
\end{equation}  
\end{corollary}

The computation of the dGSIs or dSIs of $\bo{\X}_u$ for all $u \subseteq \{1, \ldots, d\}$ using the above estimators  will require $R_{\min}$ equivalent representations of $\M$. When we are only interested in $u \subseteq \{1, \ldots, d\}$ with $|u|=1$, $R_1= \max \left(d_2,\ldots,\, d_K \right)$ equivalent representations of $\M$ are sufficient for estimating the first-order and total dGSIs or dSIs of $\X_j$ for all $j \in \{1, \ldots, d\}$.  
                 
\section{Analytical and numerical results} \label{sec:test}  
In this section, we illustrate our approach by means of three models, which allow for underlying some properties of the new indices. We derive analytical and numerical values of such indices.
   
\subsection{Linear function ($d=3$, $\NN=1$)}   
We consider
$      
    \M(\bo{\X}) = \X_1 + \X_2 +  \X_3 
$
with $\bo{\X} \sim \mathcal{N}\left(\bo{0},\, \left[\begin{array}{ccc} 
\sigma_1^2 &\rho_{12} \sigma_1\sigma_2 & \rho_{13}\sigma_1\sigma_3 \\ 
\rho_{12}\sigma_1\sigma_2 & \sigma_2^2 & \rho_{23}\sigma_2\sigma_3  \\
\rho_{13}\sigma_1\sigma_3 & \rho_{23}\sigma_2\sigma_3  & \sigma_{3}^2    
\end{array}\right]\right)$.  A dependency function of 
 $\bo{\X}$ is given by $(\X_2, \, \X_3) =r_1(\X_1, Z_2, Z_3)$ where (\cite{lamboni21,lamboni22b})
\begin{eqnarray} 
 \left\{\begin{array}{ccl}   
    X_2 &=& \frac{\rho_{12}\sigma_2}{\sigma_1}X_1 + \sqrt{1-\rho_{12}^2}  Z_2   \\
		X_3 &=& \frac{\rho_{13}\sigma_3}{\sigma_1}X_1 + \frac{\sigma_3(\rho_{23}-\rho_{12}\rho_{13})}{\sigma_2 \sqrt{1-\rho_{12}^2}} Z_2 + \sqrt{\frac{1-\rho_{12}^2-\rho_{13}^2-\rho_{23}^2+2\rho_{12}\rho_{13}\rho_{23}}{1-\rho_{12}^2}} Z_3
      \end{array} \right.   \, ,   \nonumber                       
\end{eqnarray}       
 $Z_i \sim \mathcal{N}\left(0, \sigma_i^2 \right),\, i=2,\,3$, and an equivalent representation of $\M$ is  given by (\cite{lamboni21})
\begin{eqnarray}
g_1\left(X_1, \, Z_2, \, Z_3\right)  &=& 
\left(1+\frac{\rho_{12}\sigma_2}{\sigma_1} + \frac{\rho_{13}\sigma_3}{\sigma_1}  \right)X_1  + \left(\sqrt{1-\rho_{12}^2} +  \frac{\sigma_3(\rho_{23}-\rho_{12}\rho_{13})}{\sigma_2 \sqrt{1-\rho_{12}^2}}\right) Z_2 \nonumber \\
& & + \sqrt{\frac{1-\rho_{12}^2-\rho_{13}^2-\rho_{23}^2+2\rho_{12}\rho_{13}\rho_{23}}{1-\rho_{12}^2}} Z_3   \, .    \nonumber              
\end{eqnarray}            
Based on such representation, we have the following dSIs of $\X_1$ (\cite{lamboni21}) and $(\X_1,\X_2)$   
$$
dS_1= dS_{T_1} =\frac{(\sigma_1 + \rho_{12}\sigma_2 +\rho_{13}\sigma_3)^2}{\sum_{j=1}^3\sigma_j^2 +2\rho_{12} \sigma_1\sigma_2 +2\rho_{13} \sigma_1\sigma_3 + 2\rho_{23} \sigma_2\sigma_3} \, ,
$$         
$$
dS_{12}= dS_{T_{12}} =\frac{\left(1-\rho_{12}^2 \right)(\sigma_1 + \rho_{12}\sigma_2 +\rho_{13}\sigma_3)^2 + \left(\sigma_2\left(1-\rho_{12}^2 \right) + \sigma_3(\rho_{23}-\rho_{12}\rho_{13})\right)^2 }{\left(1-\rho_{12}^2 \right)\left(\sum_{j=1}^3\sigma_j^2 +2\rho_{12} \sigma_1\sigma_2 +2\rho_{13} \sigma_1\sigma_3 + 2\rho_{23} \sigma_2\sigma_3\right)} \, ,
$$     
 respectively.  
By analogy, the two other representations of $\M$ (as $R_{\min}=3$) lead to the following results. When using $g_2\left(X_2, \, Z_3, \, Z_1\right)$, we have
$$  
dS_2= dS_{T_2} =\frac{(\sigma_2 + \rho_{12}\sigma_1 +\rho_{23}\sigma_3)^2}{\sum_{j=1}^3\sigma_j^2 +2\rho_{12} \sigma_1\sigma_2 +2\rho_{13} \sigma_1\sigma_3 + 2\rho_{23} \sigma_2\sigma_3} \, , 
$$        
$$   
dS_{23}= dS_{T_{23}} =\frac{\left(1-\rho_{23}^2 \right)(\sigma_2 + \rho_{12}\sigma_1 +\rho_{23}\sigma_3)^2 + \left(\sigma_3\left(1-\rho_{23}^2 \right) + \sigma_1(\rho_{13}-\rho_{12}\rho_{23})\right)^2 }{\left(1-\rho_{23}^2 \right)\left(\sum_{j=1}^3\sigma_j^2 +2\rho_{12} \sigma_1\sigma_2 +2\rho_{13} \sigma_1\sigma_3 + 2\rho_{23} \sigma_2\sigma_3\right)} \, .
$$   
Likewise, using $g_3\left(X_3, \, Z_1, \, Z_2\right)$ we obtain
$$
dS_3= dS_{T_3} =\frac{(\sigma_3 + \rho_{13}\sigma_1 +\rho_{23}\sigma_2)^2}{\sum_{j=1}^3\sigma_j^2 +2\rho_{12} \sigma_1\sigma_2 +2\rho_{13} \sigma_1\sigma_3 + 2\rho_{23} \sigma_2\sigma_3} \, ,
$$      
$$  
dS_{13}= dS_{T_{13}} =\frac{\left(1-\rho_{13}^2 \right)(\sigma_3 + \rho_{13}\sigma_1 +\rho_{23}\sigma_2)^2 + \left(\sigma_1\left(1-\rho_{13}^2 \right) + \sigma_2(\rho_{12}-\rho_{13}\rho_{23})\right)^2 }{\left(1-\rho_{13}^2 \right)\left(\sum_{j=1}^3\sigma_j^2 +2\rho_{12} \sigma_1\sigma_2 +2\rho_{13} \sigma_1\sigma_3 + 2\rho_{23} \sigma_2\sigma_3\right)} \, . 
$$            
        
\subsection{Portfolio model ($d=4$, $\NN=1$)}   
We consider the model given by
$        
    \M(\bo{\X}) = \X_1 \X_2 +  \X_3\X_4  
$  
with  \\ 
$(\X_1,\, \X_2) \sim \mathcal{N}_2\left(\bo{0},\, \left[\begin{array}{cc} 
\sigma_1^2 & \rho_{12}\sigma_{1} \sigma_{2} \\ 
\rho_{12}\sigma_{1} \sigma_{2} & \sigma_2^2 \end{array}\right] \right)$; 
$(\X_3,\, \X_4) \sim \bo{t}\left(\bo{0}, \, \nu, \, \left[\begin{array}{cc} 
\sigma_3^2 & \rho_{34}\sigma_{3} \sigma_{4} \\ 
\rho_{34}\sigma_{3} \sigma_{4} & \sigma_4^2 \end{array}\right] \right)$, $\nu>4$ and $(\X_1,\, \X_2)$ is independent of $(\X_3,\, \X_4)$. We are in presence of two groups of variables ($K=2$) if the correlations $\rho_{12} \neq 0$ and $\rho_{34} \neq 0$. \\     
A dependency function of $(\X_1,\, \X_2)$ is given by $\X_2 =r_1(\X_1, Z_2) = \frac{\rho_{12}\sigma_2}{\sigma_1} \X_1 + \sqrt{1-\rho_{12}^2}  Z_2$ with $\X_1 \sim \mathcal{N}\left(0, \sigma_1^2 \right)$ and  $Z_2 \sim \mathcal{N}\left(0, \sigma_2^2 \right)$.
 Likewise, we can write $\X_4 =r_3 (\X_3, Z_4)$ where (see Lemma \ref{lem:stucop}) 
$$    
 r_3 (\X_3, Z_4) = \frac{\rho_{34}\sigma_4}{\sigma_3} \X_3 + \sqrt{\frac{(1-\rho_{34}^2)\left(\nu \sigma_3^2 +\X_3^2 \right)}{\sigma_3^2(\nu +1)}}  Z_4 \, , 
$$
with $\X_3 \sim t\left(0,\, \nu,\, \sigma_3^2 \right)$ and $Z_4 \sim t\left(0,\, \nu +1,\, \sigma_4^2 \right)$.       
Thus, a  first equivalent representation of the model is given by 
\begin{eqnarray}
& & g_{13}(X_1, Z_2, \X_3, Z_4)  =
 \frac{\rho_{12}\sigma_2}{\sigma_1} \X_1^2 + \sqrt{1-\rho_{12}^2}  Z_2 \X_1 +
\frac{\rho_{34}\sigma_4}{\sigma_3} \X_3^2 + \sqrt{\frac{(1-\rho_{34}^2)\left(\nu \sigma_3^2 +\X_3^2 \right)}{\sigma_3^2(\nu +1)}}  Z_4 \X_3 \, . \nonumber   
\end{eqnarray}                   
Using the function $g_{13}$, the dSIs of $\X_1$, $(\X_1,\X_2)$, $\X_3$, $(\X_3,\X_4)$ and $(\X_1,\X_3)$ are given by 
$$
dS_1= \frac{2\rho_{12}^2\sigma_1^2\sigma_2^2}{D}, \, 
\quad 
dS_{T_1} =dS_{12}= dS_{T_{12}} = \frac{\sigma_1^2\sigma_2^2\left(1+\rho_{12}^2\right)}{D} \, , 
$$       
$$
dS_3= \frac{\rho_{34}^2 \sigma_3^2\sigma_4^2 [6(\nu-2)^2 -\nu^2(\nu-4)]}{D (\nu-4)(\nu-2)^2} \,  ,
$$
$$
dS_{T_3} = dS_{34} = dS_{T_{34}} =
\frac{\rho_{34}^2 \sigma_3^2\sigma_4^2 \frac{6(\nu-2)^2 -\nu^2(\nu-4)}{(\nu-4)(\nu-2)^2} +
\sigma_3^2\sigma_4^2\left(1-\rho_{34}^2 \right) \frac{\nu^2(\nu-4) +6 (\nu -2)}{(\nu-1)(\nu-2)(\nu-4)}}{D} \, ,   
$$               
$$
dS_{13} = dS_{1}  + dS_{3} ,\, \quad  dS_{T_{13}} =  1 \, ,  
$$     
 respectively, where 
$$
D = \sigma_1^2\sigma_2^2 \left(1+\rho_{12}^2\right) + 
\rho_{34}^2 \sigma_3^2\sigma_4^2 \frac{6(\nu-2)^2 -\nu^2(\nu-4)}{(\nu-4)(\nu-2)^2} +
\sigma_3^2\sigma_4^2\left(1-\rho_{34}^2 \right) \frac{\nu^2(\nu-4) +6 (\nu -2)}{(\nu-1)(\nu-2)(\nu-4)} \, .   
$$   
      
Using the three remaining representations of $\M$ (i.e., $R_{\min} =4$), we obtain  
$$
dS_2=  dS_1,\, \quad dS_{T_2}= dS_{T_1} \, , 
\qquad 
dS_{23}=  dS_{13},\, \quad  dS_{T_{23}}= 1 \, , 
$$      
$$
dS_4=  dS_3,\, \quad dS_{T_4}= dS_{T_3} \, , 
\qquad 
dS_{14}=  dS_{13},\, \quad  dS_{T_{14}}= 1 \, , 
$$     
$$
dS_{24}=  dS_{13},\, \quad  dS_{T_{24}}= 1 \, .  
$$            

\subsection{Multivariate g-Sobol's function with dependent variables ($d=10$, $\NN=4, K=3$)}
\label{sec:text3}  
We consider the function given by 
$$  
\M(\bo{\x}) :=  \left[ \begin{array}{c} 
\prod_{j=1}^{d=10}\frac{ |4\, \x_j \,- \,2| \,+ \,\mathcal{A}[1,j]}{1 \,+\, \mathcal{A}[1,j]} \\ 
\prod_{j=1}^{d=10}\frac{ |4\, \x_j \,-\, 2| \,+ \,\mathcal{A}[2,j]}{1 \,+ \,\mathcal{A}[2,j]} \\ 
\prod_{j=1}^{d=10}\frac{ |4 \, \x_j \,- \,2|\, + \,\mathcal{A}[3,j]}{1 \,+\, \mathcal{A}[3,j]} \\ 
\prod_{j=1}^{d=10}\frac{ |4 \, \x_j \,- \,2| \, + \,\mathcal{A}[4,j]}{1 \,+\, \mathcal{A}[4,j]} 
\end{array}      
\right] 
  \quad \mbox{with} \quad   
\mathcal{A} =\left[ \begin{array}{cccccccccc}
10 & 10 & 10 & 10 & 10 & 10 & 10 & 10 & 10 & 10 \\
20 & 20 & 20 & 20 & 20 & 20 & 20 & 20 & 20 & 20 \\ 
50 & 50 & 50 & 50 & 50 & 50 & 50 & 50 & 50 & 50 \\ 
60 & 60 & 60 & 60 & 60 & 60 & 60 & 60 & 60 & 60      
\end{array}  \right] \, .
$$      
This function includes $10$ inputs organized into $K=3$ blocks as follows: 
\begin{itemize}
\item $(\X_j \sim \mathcal{U}(0,\, 1), j=4, \ldots, 8)$ are independent variables, that is,
$\boldsymbol{\pi}_1 =\{4, \ldots, 8 \}$; 
\item  $\boldsymbol{\pi}_2 =\{1, \ldots, 3\}$ and $(\X_1, \X_2, \X_3)$ have a Gaussian copula with $\rho_{12} =0, \rho_{13} =0.01, \rho_{23} =0.85$ as the correlation values, and $\X_j \sim \mathcal{U}(0,\, 1), \, j=1, 2, 3$;
\item  $\boldsymbol{\pi}_3 =\{9, 10\}$ where $\X_9 \sim \mathcal{U}(0,\, 1), \X_{10} \sim \mathcal{U}(0,\, 1)$ with  $\X_9 + \X_{10} \leq 1$. 
\end{itemize}    
The dependency functions of $\bo{\X}_{\boldsymbol{\pi}_3}$ are given by 
$ X_{10}^c = U_{10}(1- X_{9}^c)$ and $X_{9}^c = U_9(1- X_{10}^c)$ where $U_{10}\sim \mathcal{U}(0,\, 1)$, $U_{9} \sim \mathcal{U}(0,\, 1)$;  and $U_9, \, U_{10}$ are independent of $X_{9}^c \sim Beta(1, 2)$, $X_{10}^c \sim Beta(1, 2)$ (\cite{lamboni22b}). \\
Based on the $R_{\min}=6$ representations of $\M$ (see Appendix \ref{app:equivrep}), we have computed the dGSIs using Sobol's sequences and the sample size $M=m=10000$. Table \ref{tab:ex3} contains the estimations of dGSIs, and it comes out that $X_2, \X_3, \X_9, \X_{10}$ are the most influential inputs. For fixing $\X_1$, we need the copula-based dependency model of $\bo{\X}_{\boldsymbol{\pi}_2}$, that is, $\X_1= r_1(\X_2, \Z_3, Z_1)$. Since the block of inputs $\bo{\X}_{\boldsymbol{\pi}_1}$ is not important, we have also computed dGSIs of the pairs of variables selected out of $(\bo{\X}_{\boldsymbol{\pi}_2}, \bo{\X}_{\boldsymbol{\pi}_3})$ (see Table \ref{tab:ex3}).  As expected, the total indices are always superior to the first-order indices. 
     
\begin{table}[htbp]  
\begin{center}
\begin{tabular}{lcccclcccc}
  \hline    
	\hline   
	&  \multicolumn{2}{c}{$\widehat{dGSI^{1, M}}$}	& \multicolumn{2}{c}{$\widehat{dGSI^{2, M}}$} &	&  \multicolumn{2}{c}{$\widehat{dGSI^{1, M}}$}	& \multicolumn{2}{c}{$\widehat{dGSI^{2, M}}$} \\
 & main & total & main & total &    & first-order & total & first-order & total \\
  \hline   
X1 & 0.088 & 0.090 & 0.088 & 0.090 & X1:X2 & 0.318 & 0.374 & 0.318 & 0.374 \\ 
  X2 & 0.239 & 0.296 & 0.239 & 0.295 & X1:X3 & 0.319 & 0.374 & 0.319 & 0.374 \\ 
  X3 & 0.229 & 0.284 & 0.230 & 0.284 & X1:X9 & 0.191 & 0.218 & 0.192 & 0.218 \\ 
  X4 & 0.092 & 0.093 & 0.092 & 0.093 & X1:X10 & 0.191 & 0.217 & 0.191 & 0.217 \\ 
  X5 & 0.092 & 0.093 & 0.092 & 0.093 & X2:X3 & 0.295 & 0.300 & 0.295 & 0.299 \\ 
  X6 & 0.091 & 0.094 & 0.092 & 0.093 & X2:X9 & 0.346 & 0.428 & 0.346 & 0.428 \\ 
  X7 & 0.091 & 0.093 & 0.091 & 0.093 & X2:X10 & 0.346 & 0.427 & 0.346 & 0.427 \\ 
  X8 & 0.092 & 0.094 & 0.092 & 0.094 & X3:X9 & 0.334 & 0.411 & 0.334 & 0.411 \\ 
  X9 & 0.107 & 0.133 & 0.107 & 0.133 & X3:X10 & 0.332 & 0.410 & 0.332 & 0.409 \\ 
  X10 & 0.108 & 0.134 & 0.108 & 0.134 & X9:X10 & 0.184 & 0.187 & 0.184 & 0.187 \\ 
\hline 
\hline 
\end{tabular}    
 \end{center}
  \caption{Estimates of the first-type and second-type dGSIs.}    
\label{tab:ex3}
\end{table}

\section{Conclusion} \label{sec:con}   
We have provided new dependency functions of any $d$-dimensional random vector following  the copula-based distributions with discrete variables; empirical dependency functions of inputs of complex mathematical models that comply with the imposed constraints. Combining such dependency functions with a model of interest allows for deriving different equivalent representations of such model regarding all the subsets of the model inputs. 
We used an algorithm for selecting the relevant equivalent representations of such model that allow for determining the distribution of the model outputs conditional on every subset of the model inputs, including the effects of all the subsets of the model inputs. We have extended dependent generalized sensitivity indices and dependent sensitivity indices (\cite{lamboni21}) by providing such indices for every subset of the inputs of multivariate response models, including dynamic models and single response models. \\
    
For computing dGSIs for complex models or computer codes, consistent estimators of such indices and their asymptotic distributions are provided. Analytical and numerical results confirmed  the theoretical properties of our approach such as the derivation of conditional expectations of the model outputs using equivalent representations of that model; the first-order index of any subset of inputs is less than the total index. Moreover, it came out that the sum of the first-order indices can be greater than one. In next future, it is interesting to investigate a new approach for which the first-order indices sum up to one.  
                         

\begin{appendices}
\section{Proof of Lemma \ref{lem:gaus}} \label{app:lem:gaus} 
Consider the variable $Y_k =\Phi^{-1}\left(F_k(\X_k)\right)$ if $\X_k$ is a continuous variable and $Y_k =\Phi^{-1}\left(\tau_{F_k}(X_k,\, U_k)\right)$ otherwise. It is well-known that $Y_k$
 follows the standard normal distribution with $k=1,\, \ldots,\, d$, and $\bo{Y}=(Y_j, \, \bo{Y}_{\sim j})$ has the same copula as $\bo{\X}$ (\cite{nelsen06}), as $\Phi^{-1}\circ F_k$ (resp. $\Phi^{-1}\circ \tau_{F_k}$) is a strictly increasing transformation on the range of $\X_k$.  Therefore, $\bo{Y} \sim \mathcal{N}_{d} \left(\bo{0},\, \mathcal{R}\right)$. Knowing that the dependency function of $\bo{Y}$ is given by $(Y_j, \bo{Y}_{\sim j}) =\mathcal{L} \left[Y_j,\, \bo{\Z}^\T\right]^\T$ (see \cite{lamboni21,lamboni22b}), the result follows using the inverse transformation of the form $\X_k = F_k^{\leftarrow}\circ \Phi (Y_k) $. 
  
\section{Proof of Lemma \ref{lem:stucop}} \label{app:lem:stucop}
Using the same reasoning as in Appendix \ref{app:lem:gaus}, we can see that $\bo{Y} \sim \bo{t}_d\left(\nu, \, \bo{0},\,  \mathcal{R} \right)$ with $Y_k =T_v^{-1}(F_k(\X_k))$ for a continuous variable $\X_k$ and $Y_k =T_v^{-1}(\tau_{F_k}(X_k, \, U_k))$ for a discrete variable $X_k$. We then have to derive the dependency function of $\bo{Y}$ to obtain the result, and it is done below. As  $\bo{Y} \sim \bo{t}_d(\nu, \boldsymbol{0}, \mathcal{I}) \Longleftrightarrow \mathcal{L}\bo{Y} \sim \bo{t}_d(\nu,\, \bo{0},\, \mathcal{R})$ with $\mathcal{L}$ the Cholesky factor of $\mathcal{R}$, the result holds knowing that (\cite{lamboni22b}) 
$$   
\displaystyle    
\left[\begin{array}{c}    
Y_{j} \\   
Y_{w_1}= \sqrt{\frac{\nu + Y_{j}^2}{\nu +1}} Z_{w_1}  \\         
  \vdots  \\            
 Y_{w_{d-1}}=\sqrt{\frac{\left(\nu + Y_{j}^2\right) \prod_{k=1}^{d-2}\left(\nu + k+ (Z_{w_k})^2 \right)}{\prod_{k=1}^{d-1}(\nu +k)}} Z_{w_{d-1}} \\                 
							\end{array}\right] 
											     	\sim \bo{t}_d(\nu, \boldsymbol{0}, \mathcal{I}) \,  .
$$

\section{Proof of Proposition \ref{prop:trans}} \label{app:prop:trans}   
For Point (i), we have to show that $\left(Y_{j}, \, \bo{Y}_{\sim j}\right) \stackrel{d}{=} \left(Y_{j}, T_{\sim j} \left(r_{j}\left(T_{j}^{\leftarrow}(Y_{j}),\bo{Z} \right) \right) \right)$.  For any measurable and integrable function $h: \R^{\NN+1} \to \R$, we can write
\begin{eqnarray} 
\displaystyle     
\esp\left[ h\left(Y_{j}, T_{\sim j} \left(r_{j}\left(T_{j}^{\leftarrow}(Y_{j}),\bo{Z} \right) \right) \right)  \right]  
      &=& \esp\left[ h\left(Y_{j}, T_{\sim j} \left(r_{j}\left(\X_{j},\bo{Z} \right) \right) \right)  \right] \nonumber \\   
			&=& \esp\left[ h\left(Y_{j}, \, T_{\sim j}\left(\bo{\X}_{\sim j}\right) \right)  \right] \nonumber \\ 
				&=& \esp\left[ h\left(Y_{j}, \,\bo{Y}_{\sim j}\right)  \right] = \esp\left[ h\left(\bo{Y}\right) \right] \, ,  \nonumber  
\end{eqnarray} 
bearing in mind the the theorem of transfer. Thus, Point (i) holds.\\ 
For Point (ii), first, using Point (i), we have  $ |\bo{Y}_{\sim j} | = T_{\sim j} \left(r_{j}\left(T_{j}^{\leftarrow}(|Y_{j}|),\bo{Z}  \right) \right)$. \\
Second, we can see that $R_{w_i} |Y_{w_i}| \stackrel{\text{d}}{=} Y_{w_i}$ $i=1, \ldots, d-1$, and the result holds. 
 
\section{Proof of Lemma \ref{lem:equivdist}} \label{app:lem:equivdist}
Consider any measurable and integrable function $h : \R^\NN \to \R$. As $\bo{w}_{\sim 1,k}=(\rr_k, \boldsymbol{\varpi}_k)$, it is known from Proposition \ref{prop:gdm1} that 
$$ 
\bo{\X}_{\bo{w}_{\sim 1,k}} \stackrel{d}{=} r_{w_{1,k}}\left(\X_{w_{1,k}}, \bo{Z}_{\rr_k}, \bo{Z}_{\boldsymbol{\varpi}_k}\right) \; \Longrightarrow\; \bo{\X}_{\sim \vv_k} \stackrel{d}{=} r_{\vv_k} \left(\bo{X}_{\vv_k}, \bo{Z}_{\boldsymbol{\varpi}_k} \right) \, .   
$$       
Using    
$\bo{R} := \left(\bo{\X}_{\Ori \setminus \vv_1}, \bo{\X}_{\Select \setminus \{w_{1,k_1}, \ldots, w_{1,k_m}\}},   
\bo{Z}_{\boldsymbol{\varpi}_{k_1}}, \ldots, \bo{Z}_{\boldsymbol{\varpi}_{k_m}},
\bo{Z}_{\bo{w}_{\sim 1} \setminus \{\bo{w}_{\sim 1,k_1}, \ldots, \bo{w}_{\sim 1,k_m}\}}\right) 
$ 
and the fact that the components of $(\bo{\X}_{\Ori}, \bo{\X}_{\Select}, \bo{Z}_{\bo{w}_{\sim 1}})$ are independent, we can write 
\begin{eqnarray}               
&& \esp \left[h\left(      
 g\left(\bo{\X}_{\Ori}, \bo{\X}_{\Select}, \bo{Z}_{\bo{w}_{\sim 1}} \right)  \right) | \, \bo{\X}_{\vv_1}, \X_{w_{1,k_1}}, \bo{Z}_{\rr_{k_1}}, \ldots, \X_{w_{1,k_m}}, \bo{Z}_{\rr_{k_m}} \right] \nonumber \\     
&=&  \esp_{\bo{R}} \left[h\left(   
\M\left(\bo{\X}_{\Ori}, \X_{w_{1,2}}, r_{w_{1,2}}\left(\X_{w_{1,2}}, \bo{Z}_{\bo{w}_{\sim 1,2}}\right), \ldots, \X_{w_{1,K}}, r_{w_{1,K}}\left(\X_{w_{1,K}}, \bo{Z}_{\bo{w}_{\sim 1,K}}\right)  \right)  
  \right) \right] \nonumber \\     
&=&  \esp_{\bo{R}} \left[h\left(  
\M\left(\bo{\X}_{\Ori}, \X_{w_{1,k}}, r_{w_{1,k}}\left(\X_{w_{1,k}}, \bo{Z}_{\rr_k}, \bo{Z}_{\boldsymbol{\varpi}_k}\right),\, k=2, \ldots, K \right)  \right) \right] \nonumber  \\ 
 & \stackrel{d}{=} &   \esp_{\bo{R}} \left[h\left(   
\M\left(\bo{\X}_{\Ori}, \, \bo{\X}_{\vv_k}, r_{\vv_k} \left(\bo{X}_{\vv_k}, \bo{Z}_{\boldsymbol{\varpi}_k} \right),\, k=2,\, \ldots,\, K \right)   \right)  \right] \nonumber \\  
	&= &    \esp\left[h\left(   
\M\left(\bo{\X}\right)   \right)  |\,  \bo{\X}_{\vv_1}, \bo{\X}_{\vv_{k_1}}, \ldots, \bo{\X}_{\vv_{k_m}} \right] \nonumber  \, .   
\end{eqnarray}  

\section{Proof of Lemma \ref{lem:recsets}} \label{app:lem:recsets}
For  Equation (\ref{eq:supsetpart}), at the end of the first step (i.e., $e_0=1$), $\mathcal{B}_k$ contains super-sets of $\{j_k\}$, that is, $(j_k,\, \rr_k)$ for all $j_k \in\{1, \ldots, d_k\}$ and $\rr_k \subseteq \{1, \ldots, d \}\setminus \{j_k\}$. For two super-sets $(j_{k_1}, \rr_{k_1})$ and $(j_{k_2}, \rr_{k_2})$ of $\mathcal{B}_k$, we have 
$$ 
\{j_{k_1},\, v_{1,k_1},\ldots, v_{\jmath, k_1}\} \neq \{j_{k_2},\, v_{1,k_2},\ldots, v_{\jmath, k_2}\}\; \mbox{for all} \; \jmath \in \{0 \ldots, j_{0, k}\} \, ,  \quad \mbox{and} \, 
$$ 
$$
\{v_{\jmath,k_1},\ldots, v_{d_k, k_1}\} \neq \{v_{\jmath,k_2},\ldots, v_{d_k, k_2}\} \; \mbox{for all} \; \jmath \in \{j_{0, k}+1 \ldots, d_k\} \, .     
$$
 Thus, $\left\{ u \subseteq \{1, \ldots, d_k\}  : |u|=1 \right\} \subseteq \mathcal{B}_k$. \\
Second, when $e_0 =2$ (from iteration $d_k+1$ to $\frac{d_k(d_k-1)}{2}$), we add the super-sets of $\{j_{k_1},\, j_{k_2}\}$ of the form $(j_{k_1}, j_{k_2}, \rr_{k_1,k_2})$, which were not in $\mathcal{B}_k$ at the end of the first step  $e_0=1$. As for new two super-sets, that is, $\{j_{k_1},\, j_{k_2}, \rr_{k_1,k_2}\}$, $\{j_{k_3},\, j_{k_4}, \rr_{k_3,k_4}\}$, we have $\{j_{k_1},\, j_{k_2}\} \neq \{j_{k_3},\, j_{k_4}\}$, $\{j_{k_1},\, j_{k_2}\} \neq \{j_{k_1},\, v_{1,k_1}\}$ and $ \{j_{k_1},\, v_{1,k_1}\} \neq \{j_{k_3},\, j_{k_4}\}$, the first two steps allow for obtaining $\left\{ u \subseteq \{1, \ldots, d_k\}  : 1\leq |u| \leq 2 \right\} \subseteq \mathcal{B}_k$. \\       
 Third, we repeat that procedure up to $e_0= j_{0, k}-1$ to obtain the super-sets of $\{j_{k_1},\ldots, j_{j_{0, k}-1}\}$ and 
 $\left\{ u \subseteq \{1, \ldots, d_k\}  : 1 \leq |u| \leq j_{0, k}-1 \right\} \subseteq~\mathcal{B}_k$. These operations are  possible because $\binom{d_k}{e_0} \leq \binom{d_k}{j_{0, k}}$ for all $e_0=1,\ldots, j_{0, k}-1$, and we avoid permutations ($\bo{w}_k$) that bring replicated sets in both $\mathcal{B}_k$ and $\mathcal{E}_k$. \\  
Fourth, the iterations $\binom{d_k}{j_{0, k}-1} < i \leq \binom{d_k}{j_{0, k}}$ (when possible) aim to add the remaining subsets of $j_{0, k}$ elements. \\  
Fifth, we have $\left\{ u \subseteq \{1, \ldots, d_k\}  : |u|=j_{0, k} +1 \right\} \subseteq \mathcal{B}_k$ because for any $\vv_1 \subseteq \{1, \ldots, d_k\}$ with $|\vv_1|=j_{0, k} +1$,  there exists $\bo{w}_k^* \in \mathcal{P}_k$ such that $\{w_{_{j_{0, k}+2,k}}^*, \ldots, w_{_{d_k, k}}^* \} \bigcap \vv_1=\emptyset$. Indeed, $\bo{w}_k^*$ was added in $\mathcal{P}_k$ when constructing all the subsets $u \subseteq \{1, \ldots, d_k\}$ with $|u|= d_k-j_{0, k}-1 < j_{0, k}$ thanks to  $\mathcal{E}_k$ and the fact that $\binom{d_k}{j_{0, k}+1} =  \binom{d_k}{d_k-j_{0, k}-1}$. Thus, $\vv_1 =\{w_{_{1,k}}^*, \ldots, w_{_{j_{0, k} +1, k}}^*\}$.    
Finally, we use the same reasoning to obtain the results. \\ 
Equation (\ref{eq:supsetpart1}). is obvious by construction (see Algorithm \ref{algo:sets}).
    
\section{Proof of Theorem \ref{theo:tierep}} \label{app:theo:tierep}
For Point (i), first, for $\vv_k \subseteq \{1, \ldots, d_k\}$ with $|\vv_k|>0$, there exists $\bo{w}_k^* \in \mathcal{P}_k$ such that $\vv_k =\{w_{1, k}^*,  \ldots, w_{|\vv_k|, k}^* \}$ according to Lemma \ref{lem:recsets}.    
Lemma \ref{lem:equivdist} ensures the  determination of the distribution of $\M(\bo{\X})$ conditional on $\bo{\X}_{\vv_k}$ using $g$ associated with $\bo{w}_k^*$. \\  
Second,  for $u := (\vv_1, \vv_2, \ldots, \vv_K)$ where $\vv_1 \subseteq \Ori,\; \vv_k \subseteq \{1, \ldots, d_k\}$ and $|\vv_k| >0$ with $k=2, \ldots, K$, there exists only one permutation $\bo{w}_k^* \in \mathcal{P}_k$ such that $\vv_k =\{w_{1, k}^*,  \ldots, w_{|\vv_k|, k}^* \},\, \forall\, k=2, \ldots, K$. 
As only one representation of $\M$ associated with $\bo{w}_k^*,\, k=2,\ldots, K$ allows for determining the distribution of $\M(\bo{\X})$ conditional on $\bo{\X}_{u}$, and  $|\mathcal{P}_k| = \binom{d_k}{j_{0, k}}$, then $R := \prod_{ k =2}^K  \binom{d_k}{j_{0, k}}$ different representations of $\M$ are needed to obtain the distribution of $\M(\bo{\X})$ conditional on $\bo{\X}_u$ for all $u \subseteq \{1, \ldots, d\}$. The result follows because $R$ is the highest number of possibilities of  \\
$\left\{\vv_k \subseteq \{1, \ldots, d_k\},\, k=2, \ldots, K : |\vv_k| =j_{0, k} \right\}$, and other possibilities are in the $R$ representations (see Lemma \ref{lem:recsets}). \\ 
For Point (ii),  let $R_{p_2\ldots p_K} =\max_{2\leq k \leq K} \binom{d_{k}}{p_{k}}$. Using Lemmas \ref{lem:equivdist}-\ref{lem:recsets}, $\binom{d_k}{p_k}$ equivalent representations of $\M$ are necessary to determine the distribution of $\M(\bo{\X})$ conditional on $\bo{\X}_{\vv_{k}}$ for all $\vv_{k} \subseteq \{1,\ldots, d_{k}\}$ with $|u_{k}| \leq p_{k}$. With $R_{p_2\ldots p_K}$ representations, we can assess the effect of any groups of variables because $\binom{d_k}{p_k} \leq R_{p_2\ldots p_K}, \, \forall\,  k=2,\ldots, K$. 
	       
\section{Proof of Proposition \ref{prop:gsipro}} \label{app:prop:gsipro}
The proofs are straightforward because each equivalent representation of $\M$ includes only independent variables. The results rely on the Hoeffding decomposition of such  equivalent representation. The proofs are similar to those provided in \cite{lamboni21} (see Proposition 2).
                     
\section{Proof of Theorem \ref{theo:sigestM}} \label{app:theo:sigestM}
The proofs are straightforward because each equivalent representation of $\M$ includes only independent variables. By considering such equivalent representation and knowing that the kernels used in this paper are particular cases of those provided in \cite{lamboni18a} (see Theorem 1), we obtain MVUEs bearing in mind the U-statistic theory.  Indeed, each kernel remains unchanged when one permutes $\bo{\X}_{i, \oo\cap u}^{(1)}$ with $\bo{\X}_{i, \oo\cap u}^{(2)}$ or  $\bo{\X}_{i, \oo\setminus u}^{(1)}$ with $\bo{\X}_{i, \oo\setminus u}^{(2)}$ or $\bo{\X}_{i, \Select \cap u}^{(1)}$ with $\bo{\X}_{i, \Select \cap u}^{(2)}$ or $\bo{\X}_{i, \Select \setminus u}^{(1)}$ with $\bo{\X}_{i, \Select \setminus u}^{(2)}$ or $\bo{Z}_{i, \rr}^{(1)}$ with $\bo{Z}_{i, \rr}^{(2)}$ or
 $\bo{Z}_{i, \boldsymbol{\varpi}}^{(1)}$ with $\bo{Z}_{i, \boldsymbol{\varpi}}^{(2)}$. More details can be found in \cite{lamboni18a} (see Theorems 1-3) and \cite{lamboni18}.
        
\section{Proof of Theorem \ref{theo:estgsiM}} \label{app:theo:estgsiM}
The results about the consistency are obtained by applying the Slutsky theorem.\\  
For the asymptotic distributions of Point (i), the derivation of these results are similar to those of Theorem 4 provided in \cite{lamboni18a}. The proofs of the asymptotic distributions of Point (ii) are similar to the proofs of Theorem 6 provided in \cite{lamboni22} under the condition $m/M \to 0$. 

\section{Equivalent representations of the function used in Section~\ref{sec:text3}} \label{app:equivrep}
The $R_{\min} =6$ equivalent representations of $\M(\bo{\X})$ are given as follows:
\begin{eqnarray} 
g_1 := \M\left(\X_1, r_1(\X_1, Z_2, Z_3), \bo{\X}_{\boldsymbol{\pi}_1}, \, \X_9, r_9(\X_9, Z_{10}) \right) & \mbox{for} & \X_1, \X_9, (\X_1, \X_2), (\X_1, \X_9),  \nonumber \\
& & \X_k, (\X_1, \X_k), \forall \, k \in \boldsymbol{\pi}_1, \ldots \, ; \nonumber \\
g_2 := \M\left(\X_2, r_2(\X_2, Z_3, Z_1), \bo{\X}_{\boldsymbol{\pi}_1}, \, \X_9, r_9(\X_9, Z_{10}) \right)  &\mbox{for}  & \X_2, (\X_2, \X_3), (\X_2, \X_9), (\X_9, \X_{10}), \ldots \, ; \nonumber \\
g_3 := \M\left(\X_3, r_3(\X_3, Z_1, Z_2), \bo{\X}_{\boldsymbol{\pi}_1}, \, \X_9, r_9(\X_9, Z_{10}) \right) &\mbox{for}  & \X_3, (\X_3, \X_1), (\X_3, \X_9), \ldots \, ; \nonumber \\
g_4 := \M\left(\X_1, r_1(\X_1, Z_2, Z_3), \bo{\X}_{\boldsymbol{\pi}_1}, \, \X_{10}, r_{10}(\X_{10}, Z_{9}) \right) & \mbox{for}  &  \X_{10}, (\X_1, \X_{10}), \ldots \, ; \nonumber \\
g_5 := \M\left(\X_2, r_2(\X_2, Z_3, Z_1), \bo{\X}_{\boldsymbol{\pi}_1}, \, \X_{10}, r_{10}(\X_{10}, Z_{9}) \right) &\mbox{for} &  (\X_2, \X_{10}), \ldots \, ; \nonumber \\
g_6 := \M\left(\X_3, r_3(\X_3, Z_1, Z_2), \bo{\X}_{\boldsymbol{\pi}_1}, \, \X_{10}, r_{10}(\X_{10}, Z_{9}) \right) &\mbox{for} & (\X_3, \X_{10}), \ldots \, , \nonumber 
\end{eqnarray}        
where $r_1, r_2, r_3$ are particular cases of DMs provided in Lemma \ref{lem:gaus}.
         
\end{appendices}  
                       
\section*{References}         
\bibliographystyle{elsarticle-num}        

\begin{thebibliography}{10}
\expandafter\ifx\csname url\endcsname\relax
  \def\url#1{\texttt{#1}}\fi
\expandafter\ifx\csname urlprefix\endcsname\relax\def\urlprefix{URL }\fi
\expandafter\ifx\csname href\endcsname\relax
  \def\href#1#2{#2} \def\path#1{#1}\fi

\bibitem{sobol93}
I.~M. Sobol, Sensitivity analysis for non-linear mathematical models,
  Mathematical Modelling and Computational Experiments 1 (1993) 407--414.

\bibitem{sobol01}
I.~M. Sobol, Global sensitivity indices for nonlinear mathematical models and
  their {Monte Carlo} estimates, Mathematics and Computers in Simulation 55
  (2001) 271 -- 280.

\bibitem{lamboni18a}
M.~Lamboni, Multivariate sensitivity analysis: Minimum variance unbiased
  estimators of the first-order and total-effect covariance matrices,
  Reliability Engineering \& System Safety 187 (2019) 67 -- 92.

\bibitem{lamboni11}
M.~Lamboni, H.~Monod, D.~Makowski, Multivariate sensitivity analysis to measure
  global contribution of input factors in dynamic models, Reliability
  Engineering and System Safety 96 (2011) 450--459.

\bibitem{gamboa14}
F.~Gamboa, A.~Janon, T.~Klein, A.~Lagnoux, Sensitivity indices for multivariate
  outputs, Comptes Rendus Mathematique 351~(7) (2013) 307--310.

\bibitem{gamboa14b}
F.~Gamboa, A.~Janon, T.~Klein, A.~Lagnoux, Sensitivity analysis for
  multidimensional and functional outputs, Electron. J. Statist. 8~(1) (2014)
  575--603.

\bibitem{xiao17}
S.~Xiao, Z.~Lu, L.~Xu, Multivariate sensitivity analysis based on the direction
  of eigen space through principal component analysis, Reliability Engineering
  \& System Safety 165 (2017) 1 -- 10.

\bibitem{lamboni19}
M.~Lamboni, Derivative-based generalized sensitivity indices and {S}obol'
  indices, Mathematics and Computers in Simulation 170 (2020) 236 -- 256.

\bibitem{perrin21}
T.~Perrin, O.~Roustant, J.~Rohmer, O.~Alata, J.~Naulin, D.~Idier, R.~Pedreros,
  D.~Moncoulon, P.~Tinard, Functional principal component analysis for global
  sensitivity analysis of model with spatial output, Reliability Engineering \&
  System Safety 211 (2021) 107522.

\bibitem{daveiga09}
S.~D. Veiga, F.~Wahl, F.~Gamboa, Local polynomial estimation for sensitivity
  analysis on models with correlated inputs, Technometrics 51 (2009) 452--463.

\bibitem{mara12}
T.~A. Mara, S.~Tarantola, Variance-based sensitivity indices for models with
  dependent inputs, Reliability Engineering \& System Safety 107 (2012) 115 --
  121.

\bibitem{kucherenko12}
S.~Kucherenko, S.~Tarantola, P.~Annoni, Estimation of global sensitivity
  indices for models with dependent variables, Computer Physics Communications
  183~(4) (2012) 937 -- 946.

\bibitem{hao13}
W.~Hao, L.~Zhenzhou, P.~Wei, Uncertainty importance measure for models with
  correlated normal variables, Reliability Engineering \& System Safety 112
  (2013) 48 -- 58.

\bibitem{chastaing12}
G.~Chastaing, F.~Gamboa, C.~Prieur, Generalized {H}oeffding-{S}obol'
  decomposition for dependent variables - applications to sensitivity analysis,
  Electronic Journal of Statistics 6 (2012) 2420--2448.

\bibitem{kucherenko17}
S.~Kucherenko, O.~Klymenko, N.~Shah, {S}obol' indices for problems defined in
  non-rectangular domains, Reliability Engineering \& System Safety 167 (2017)
  218 -- 231.

\bibitem{mara15}
T.~A. Mara, S.~Tarantola, P.~Annoni, Non-parametric methods for global
  sensitivity analysis of model output with dependent inputs, Environmental
  Modelling \& Software 72 (2015) 173 -- 183.

\bibitem{tarantola17}
S.~Tarantola, T.~A. Mara, Variance-based sensitivity indices of computer models
  with dependent inputs: The fourier amplitude sensitivity test, International
  Journal for Uncertainty Quantification 7~(6) (2017) 511--523.

\bibitem{lamboni21}
M.~Lamboni, S.~Kucherenko, Multivariate sensitivity analysis and
  derivative-based global sensitivity measures with dependent variables,
  Reliability Engineering \& System Safety 212 (2021) 107519.

\bibitem{skorohod76}
A.~V. Skorohod, On a representation of random variables, Theory Probab. Appl
  21~(3) (1976) 645--648.

\bibitem{lamboni22b}
M.~Lamboni, Efficient dependency models: simulating dependent random variables,
  Mathematics and Computers in Simulation , submitted on 03/01/2021.

\bibitem{owen14b}
A.~Owen, Sobol' indices and shapley value, Journal on Uncertainty
  Quantification 2 (2014) 245--251.

\bibitem{ferguson67}
T.~Ferguson, Mathematical Statistics: A Decision Theoretic Approach, Academic
  Press, New York, 1967.

\bibitem{ruschendorf81}
L.~Rüschendorf, Stochastically ordered distributions and monotonicity of the
  oc-function of sequential probability ratio tests, Series Statistics 12~(3)
  (1981) 327--338.

\bibitem{ruschendorf05}
L.~Rüschendorf, Stochastic ordering of risks, influence of dependence and a.s.
  constructions, in: Advances on Models, Characterizations and Applications, N.
  Balakrishnan, I. G. Bairamov, O. L. Gebizlioglu, 2005.

\bibitem{ruschendorf09}
L.~Rüschendorf, On the distributional transform, sklar's theorem, and the
  empirical copula process, Journal of Statistical Planning and Inference
  139~(11) (2009) 3921 -- 3927.

\bibitem{rosenblatt52}
M.~Rosenblatt, Remarks on a multivariate transformation, Ann. Math. Statist.
  23~(3) (1952) 470--472.

\bibitem{obrien75}
G.~L. O'Brien, The comparison method for stochastic processes, The Annals of
  Probability 3~(1) (1975) 80 -- 88.

\bibitem{arjas78}
E.~Arjas, T.~Lehtonen, Approximating many server queues by means of single
  server queues, Mathematics of Operations Research 3 (1978) 205--223.

\bibitem{ruschendorf93}
L.~Rüschendorf, V.~{de Valk}, On regression representations of stochastic
  processes, Stochastic Processes and their Applications 46~(2) (1993) 183 --
  198.

\bibitem{nelsen06}
R.~Nelsen, An Introduction to Copulas, Springer-Verlag, New York., 2006.

\bibitem{mcneil15}
A.~J. McNeil, R.~Frey, P.~Embrechts, Quantitative Risk Management, Princeton
  University Press, Princeton and Oxford, 2015.

\bibitem{rosenblatt56}
M.~Rosenblatt, Remarks on some nonparametric estimates of a density function,
  Annals of Mathematical Statistics 27 (1956) 832--837.

\bibitem{parzen62}
E.~Parzen, On estimation of a probability density function and mode, Annals of
  Mathematical Statistics 33 (1962) 1065--1076.

\bibitem{epanechnikov69}
V.~Epanechnikov, Nonparametric estimation of a multidimensional probability
  density, Theory Probab. Appl. 14 (1969) 153--158.

\bibitem{silverman86}
B.~Silverman, Density Estimation for Statistics and Data Analysis, Chapman \&
  Hall, New York, 1986.

\bibitem{clayton78}
D.~G. Clayton, A model for association in bivariate life tables and its
  application in epidemiological studies of familial tendency in chronic
  disease incidence, Biometrika 65~(1) (1978) 141--151.

\bibitem{joe97}
H.~Joe, Multivariate models and dependence concepts, Boca Raton ; London; New
  York: Chapman \& Hall/CRC., 1997.

\bibitem{smith10}
M.~Smith, A.~Min, C.~Almeida, C.~Czado, Modeling longitudinal data using a
  pair-copula decomposition of serial dependence, Journal of the American
  Statistical Association 105~(492) (2010) 1467--1479.

\bibitem{durante15}
F.~Durante, C.~Sempi, Principles of copula theory, CRC/Chapman \& Hall, London,
  2015.

\bibitem{koenker78}
R.~Koenker, G.~Bassett, Regression quantiles, Econometrica 46 (1978) 33--50.

\bibitem{truong89}
Y.~K. Truong, Asymptotic properties of kernel estimators based on local
  medians, The Annals of Statistics 17~(2) (1989) 606--617.

\bibitem{hendricks92}
W.~Hendricks, R.~Koenker, Hierarchical spline models for conditional quantiles
  and the demand for electricity, Journal of the American Statistical
  Association 87~(417) (1992) 58--68.

\bibitem{koenker94}
R.~Koenker, P.~Ng, S.~Portnoy, {Quantile smoothing splines}, Biometrika 81~(4)
  (1994) 673--680.

\bibitem{koenker01}
R.~Koenker, K.~F. Hallock, Quantile regression, The Journal of Economic
  Perspectives 15~(4) (2001) 143--156.

\bibitem{bassett82}
G.~B. Jr., R.~Koenker, An empirical quantile function for linear models with
  iid errors, Journal of the American Statistical Association 77~(378) (1982)
  407--415.

\bibitem{koenker05}
R.~Koenker, Quantile Regression, Cambridge University Press, 2005.

\bibitem{takeuchi06}
I.~Takeuchi, Q.~V. Le, T.~D. Sears, A.~J. Smola, Nonparametric quantile
  estimation, Journal of Machine Learning Research 7~(45) (2006) 1231--1264.

\bibitem{currin91}
C.~Currin, T.~Mitchell, M.~Morris, D.~Ylvisaker, Bayesian prediction of
  deterministic functions, with applications to the design and analysis of
  computer experiments, Journal of the American Statistical Association
  86~(416) (1991) 953--963.

\bibitem{haylock96}
R.~G. Haylock, A.~O\^Hagan, J.~M. Bernardo, On inference for outputs of
  computationally expensive algorithms with uncertainty on the inputs,
  international meeting; 5th, bayesian statistics, in: Bayesian statistics,
  International meeting; 5th, Bayesian statistics, Vol.~5, 1996, pp. 629--638.

\bibitem{kennedy01}
M.~C. Kennedy, A.~O'Hagan, Bayesian calibration of computer models, Journal of
  the Royal Statistical Society: Series B (Statistical Methodology) 63~(3)
  (2001) 425--464.

\bibitem{jeremy04}
J.~E. Oakley, A.~O'Hagan, Probabilistic sensitivity analysis of complex models:
  a bayesian approach, Journal of the Royal Statistical Society: Series B
  (Statistical Methodology) 66~(3) (2004) 751--769.

\bibitem{lamboni18}
M.~Lamboni, Uncertainty quantification: a minimum variance unbiased (joint)
  estimator of the non-normalized {S}obol' indices, Statistical Papers 61
  (2020) 1939--1970.

\bibitem{lamboni22}
M.~Lamboni, Weak derivative-based expansion of functions: Anova and some
  inequalities, Mathematics and Computers in Simulation 194 (2022) 691--718.

\bibitem{lamboni16}
M.~Lamboni, Global sensitivity analysis: an efficient numerical method for
  approximating the total sensitivity index, International Journal for
  Uncertainty Quantification 6~(1) (2016) 1--17.

\bibitem{lamboni16b}
M.~Lamboni, Global sensitivity analysis: a generalized, unbiased and optimal
  estimator of total-effect variance, Statistical Papers 59~(1) (2018)
  361--386.

\bibitem{lamboni20}
M.~Lamboni, Derivative-based integral equalities and inequality: A
  proxy-measure for sensitivity analysis, Mathematics and Computers in
  Simulation 179 (2021) 137 -- 161.

\end{thebibliography}
                         

\end{document}